\documentclass{amsart}
\usepackage{amsmath,amssymb}
\usepackage{xy}
\xyoption{all}
\newcommand{\ie}{{\em i.e., }}
\newcommand{\eg}{{\em e.g., }}
\newcommand{\A}{\ensuremath{\mathcal A}}

\newcommand{\C}{\ensuremath{\mathbb C}}
\newcommand{\CC}{\ensuremath{\mathcal C}}

\newcommand{\E}{\ensuremath{\mathcal E}}
\newcommand{\F}{\ensuremath{\mathcal F}}

\newcommand{\K}{\ensuremath{\mathcal K}}

\renewcommand{\O}{\ensuremath{\mathcal O}}

\newcommand{\R}{\ensuremath{\mathbb R}}
\renewcommand{\S}{\ensuremath{\mathcal S}}
\newcommand{\SC}{\ensuremath{\mathcal S \mathcal C}}

\newcommand{\Z}{\ensuremath{\mathbb Z}}

\newcommand{\Cliff}{{\operatorname{Cliff}}}

\newcommand{\cs}{{C^{*}}}

\newcommand{\into}{\hookrightarrow}

\renewcommand{\epsilon}{{\varepsilon}}

\newcommand{\gtimes}{{\hat{\otimes} \,}}
\newcommand{\Hom}{\operatorname{Hom}}
\newcommand{\ideal}{\vartriangleleft}
\newcommand{\id}{\operatorname{id}}

\newcommand{\spinc}{{spin}$^c$\!}
\newcommand{\Spin}{\operatorname{{Spin}}}

\newcommand\<{\langle}
\renewcommand\>{\rangle}

\theoremstyle{plain}
\newtheorem{thm}{Theorem}[section] 

\newtheorem{cor}[thm]{Corollary}

\newtheorem{lem}[thm]{Lemma}

\newcommand{\Pf}{{\em Proof}. }

\newtheorem{prop}[thm]{Proposition}

\newtheorem{thm*}{Theorem}
\newcommand{\EPf}{\hfill $\Box$}
  {\addvspace{\bigskipamount}\noindent{\bf Example\ }}%
  {}

\newtheorem{ThomThm}[thm]{Thom Isomorphism Theorem}

{\begin{list}{$\bullet$}{\parsep=0pt \itemsep=0pt 
                         \topsep=6pt \partopsep=\parskip}}%
{\end{list}}

\newcounter{ictr}

\newcounter{nctr}

\theoremstyle{definition}
\newtheorem{dfn}[thm]{Definition}
\newtheorem{example}[thm]{Example}
\newtheorem{examples}[thm]{Examples}

\title[On $C^*$-Algebras and $K$-theory for Fredholm Manifolds]
{On $C^*$-Algebras and $K$-theory for Infinite-Dimensional Fredholm
Manifolds}\author[Dorin Dumitra\c scu and Jody Trout]{Dorin Dumitra\c scu and Jody Trout$^\dag$}
\address{6188 Bradley Hall\\ 
         Dartmouth College\\
         Hanover, NH 03755}
\address{617 N. Santa Rita\\
University of Arizona\\
Tucson, AZ 85721-0089}
\thanks{$^\dag$ The second author was partially supported by NSF Grant DMS-0071120}

\subjclass{19, 46, 47, 55, 57, 58}
\keywords{$C^*$-algebra, Fredholm manifold, direct limit, $K$-theory, $K$-homology,  Poincar\'e duality}
\begin{document}

\begin{abstract}
Let $M$ be a smooth Fredholm manifold modeled on a separable infinite-dimensional Euclidean space $\E$ with Riemannian metric $g$. Given an augmented Fredholm filtration $\F $
 of  $M$ by finite-dimensional submanifolds $\{M_n\}_{n=k}^\infty$, we associate to the triple $(M, g, \F)$ a non-commutative direct limit $\cs$-algebra
\begin{equation*}
\A(M, g, \F) = \varinjlim \, \A(M_n)
\end{equation*} 
that can play the role of the algebra of functions vanishing at infinity on the non-locally compact space $M$. The $\cs$-algebra $\A(\E)$, as constructed by Higson-Kasparov-Trout  for their Bott periodicity theorem, is isomorphic to our construction when $M = \E$. If $M$ has an oriented $\Spin_q$-structure $(1 \leq q \leq \infty)$, then the $K$-theory of this $\cs$-algebra is the same (with dimension shift) as the topological $K$-theory of $M$ defined by Mukherjea. Furthermore, there is a Poincar\'e duality isomorphism of this $K$-theory of $M$ with the compactly supported $K$-homology of $M$, just as in the finite-dimensional spin setting. 
\end{abstract}

\maketitle

 \section{Introduction}
 
Infinite-dimensional Hilbert manifolds have been studied since the 1960's, with main applications in infinite-dimensional differential topology, global analysis, non-linear PDEs, and other areas. This paper is concerned with constructing $\cs$-algebras and computing the $K$-theory for a particular class of infinite-dimensional Hilbert manifolds, namely {\it Fredholm manifolds} \cite{Ee,EeElw1,ElwTr}. This is part of a research program to introduce concepts and techniques from Alain Connes' noncommutative geometry \cite{Con94}, \eg $\cs$-algebras, $K$-theory, cyclic (co)homology, and spectral triples, into the study of Fredholm manifolds.
  
But first, let us review the finite-dimensional case. Given $M$ a {\it finite-dimensional} Riemannian manifold, let $C_0(M)$ be the commutative $\cs$-algebra of all continuous complex-valued functions which vanish at infinity on $M$. This $\cs$-algebra categorically encodes the topological properties of $M$ \cite{WO93} and, by the Serre-Swan theorem,  plays a dual role in the $K$-theory of $M$:
\begin{equation*}
K^j(M) \cong K_j(C_0(M)) , \quad j = 0, 1,
\end{equation*}
where $K^j(M)$ is the (reduced) topological $K$-theory of $M$ \cite{Ati67}. Furthermore, if $M$ has a spin (or spin$^c$) structure \cite{LM89}, there is a Poincar\'e duality isomorphism \cite{HR00,Ren03}:
\begin{equation*}
K^{n-j}(M) \cong K_j^c(M), \quad j = 0, 1,
\end{equation*}
where $K_j^c(M)$ denotes the dual (compactly supported) $K$-homology of $M$ and $n$ is the dimension of $M$. 

The other $\cs$-algebra for a finite-dimensional $M$ is non-commutative and constructed using the Riemannian metric $g$. For each $x \in M$, the tangent space $T_xM$ of $M$ is a finite-dimensional Euclidean space with inner product $g_x$. Thus, we can form the complex Clifford algebra $\Cliff(T_xM, g_x)$ (see Section 2). It has a canonical structure as a finite-dimensional $\Z_2$-graded $\cs$-algebra. The family of $\cs$-algebras $\{\Cliff(T_xM, g_x)\}_{x \in M}$ naturally forms a $\Z_2$-graded, $\cs$-algebra vector bundle $\Cliff(TM) \to M$, called the {\it Clifford algebra bundle} of $M$ \cite{ABS64}. We then can define
\begin{equation*}
\CC(M) = C_0(M, \Cliff(TM))
\end{equation*}
to be the $\cs$-algebra of continuous sections of the Clifford algebra bundle of $M$ vanishing at infinity. This $\cs$-algebra was used by Kasparov \cite{Kas88} in studying the Novikov Conjecture, where he used the notation $\CC_\tau(M)$. If $M$ is even-dimensional and has a spin structure (or, more generally, a \spinc-structure) then this $\cs$-algebra is Morita equivalent to $C_0(M)$. (In general, $\CC(M)$ is Morita equivalent to $C_0(TM)$.) By the Morita invariance of $K$-theory, it follows that
\begin{equation*}
K_j(\CC(M)) \cong K_j(C_0(M)) \cong K^j(M), \quad j = 0, 1.
\end{equation*}
For $M$ odd-dimensional and spin, this is more complicated. (See Proposition \ref{prop:morita}.)

If $M$ is an {\it infinite-dimensional} Hilbert manifold \cite{Lang}, modeled on a separable infinite-dimensional Euclidean (\ie real Hilbert) space $\E$, then these two constructions do not  work. Both fail since compact subsets of $M = \E$ are ``thin'', \ie contained in finite-dimensional subspaces. Thus, $C_0(\E) = \{0\}$ since there are no compactly supported continuous functions on $\E$ which are non-zero. However, the Clifford $\cs$-algebra has been generalized by Higson-Kasparov-Trout \cite{HKT98} to the case $M = \E$, by a direct limit construction that exploits an important property  of Clifford algebras with respect to orthogonal sums (see equation (\ref{eqn:orthogonality})). The component $\cs$-algebras in the direct limit are given by
\begin{equation*}
\A(E^a) = C_0(\R) \gtimes  \CC(E^a) \cong C_0(\R) \gtimes  C_0(E^a, \Cliff(E^a)) 
\end{equation*}
where $\gtimes$ denote the $\Z_2$-graded tensor product \cite{Bla98} and $C_0(\R)$ is graded by even and odd functions. Since the map $E^a \mapsto \A(E^a)$ is functorial with respect to inclusions of finite-dimensional susbspaces, one can construct a non-commutative direct limit $\cs$-algebra (in the better notation of \cite{HK01}):
\begin{equation*}
\A(\E) = \varinjlim_{E^a \subset \E} \A(E^a)
\end{equation*}
where the direct limit is taken over {\it all finite-dimensional subspaces} $E^a \subset \E$. 
(See Example \ref{ex:HKT} for more on this construction and how it fits into our theory.) This $\cs$-algebra was used to prove an equivariant Bott periodicity theorem for infinite-dimensional Euclidean spaces \cite{HKT98} and  has had applications to proving cases of the Novikov Conjecture and, more generally, the Baum-Connes Conjecture \cite{HK01,Yu00}.

Now, suppose the Hilbert manifold $M$ is fibered as the total space of a smooth infinite rank  Euclidean vector bundle $p : F \to X$, with fiber $\E$ and compatible affine connection $\nabla$, over a finite-dimensional Riemannian manifold $X$. Let $p_a : F^a \to X$ be a {\it finite rank} subbundle of $F$. Using the connection $\nabla$ and the metrics on $F$ and $X$, we can give the total space $F_a$ a canonical structure of a Riemannian manifold and define the component $\cs$-algebra
\begin{equation*}
\A(F^a) = C_0(\R) \gtimes  \CC(F^a) \cong C_0(\R) \gtimes  C_0(F^a, \Cliff(TF^a)). 
\end{equation*} 
Since the map $F^a \mapsto \A(F^a)$ is functorial with respect to inclusions of finite-dimensional subbundles \cite{Trou03}, we can then construct a direct limit $\cs$-algebra:
\begin{equation*}
\A(F, \nabla) = \varinjlim_{F^a \subset F} \A(F^a)
\end{equation*}
where the direct limit is taken over {\it all finite rank subbundles} $p_a : F^a \to X$ of $F$. Trout \cite{Trou03} used this $\cs$-algebra to prove an equivariant Thom isomorphism theorem for infinite rank Euclidean bundles, which reduces to the Higson-Kasparov-Trout Bott periodicity theorem when the base manifold $X$ is a point.

For a more general {\it curved} Hilbert manifold $M$, with Riemannian metric $g$,  there does not seem to be a natural generalization of the previous constructions. Based on the above, one would be tempted to construct a direct limit $\cs$-algebra
\begin{equation*}
\text{``} \A(M) = \varinjlim_{M_a \subset M} \A(M_a) \text{''}
\end{equation*}
where the component $\cs$-algebras should be given by
\begin{equation*}
\A(M_a) = C_0(\R) \gtimes  \CC(M_a)
\end{equation*}
and the direct limit is taken over {\it all finite-dimensional submanifolds} $M_a \subset M$. The problem is that, even though the component $\cs$-algebras have many functoriality properties (as discussed in Section 2), if we are given smooth (isometric) inclusions
\begin{equation*}
M_a \subset M_b \subset M_c
\end{equation*}
of finite-dimensional submanifolds of $M$, there is no obvious way to define a commuting diagram (as there is in the Bott periodicity and Thom isomorphism cases)
\begin{equation}\label{diag:functor}
\xymatrix{ & \A(M_b) \ar[dr] &  \\ \A(M_a) \ar[ur] \ar[rr]  &  & \A(M_c)}
\end{equation}
needed to construct the corresponding direct limit. 

However, if the Hilbert manifold $M$ has a {\it Fredholm structure}, then we can construct a direct limit $\cs$-algebra by choosing an appropriate {\it countable sequence} $\{M_n\}_{n=k}^\infty$ of expanding, topologically closed, finite-dimensional submanifolds of $\dim(M_n) = n$. The sequence $\{M_n\}_{n=k}^\infty$ is called a {\it Fredholm filtration} of $M$. (See Section 3 for the geometric definitions and details.) The countability of this sequence of submanifolds clearly simplifies the direct limit construction since only each ``Gysin'' map $\A(M_n) \to \A(M_{n+1})$ needs to be constructed, which will require some non-trivial geometry (\ie connections and normal bundles.)

Equip the Riemannian Fredholm manifold $(M, g)$ with an {\it augmented} Fredholm filtration  $\F = (M_n , U_n)_{n= k}^\infty$ (as in Definition \ref{D:aFf}) where $U_n$ is a total open tubular neighborhood of $M_n \hookrightarrow M_{n+1}$. Section 4 contains the construction of a noncommutative direct limit $\cs$-algebra for the triple $(M, g, \F)$:
\begin{equation*}
\A(M, g, \F) = \varinjlim\limits  \, \A(M_n)
\end{equation*} 
that can play the role of the algebra of functions vanishing at infinity on $M$. 

Using ideas of Mukherjea \cite{Mkhr1,Mkhr2} to associate cohomology functors to Fredholm manifolds via Fredholm filtrations, the {\it topological $K$-theory groups} of $(M, \F)$ are defined as the direct limit:
\begin{equation*}
K^{\infty-j}(M, \F) = \varinjlim \, K^{n-j}(M_n), \quad j = 0, 1,
\end{equation*}
where the connecting map $K^{n-j}(M_n) \to K^{(n+1) - j}(M_{n+1})$ is the Gysin (or shriek) map (Definition \ref{D:topKth}) of the embedding $M_n \hookrightarrow M_{n+1}$, and the inspiration for our connecting map $\A(M_n) \to \A(M_{n+1})$. Note that this definition does, in general, depend on the choice of Fredholm filtration, since the sequence $\{M_n\}_{n=k}^\infty$ may not be $K$-orientable \cite{Dyer69,CS84}.

But, using appropriate notions of {$\Spin_q$-structures} (see Section 5.2) for Riemannian Fredholm manifolds, originally investigated by Anastasiei \cite{Anst} and de la Harpe \cite{dlH}, the following Serre-Swan and Poincar\'e duality isomorphism theorem (combining Theorems \ref{Thm:FMPD} and \ref{T:Kth-A(M)}) is obtained:
\begin{thm} Let $(M, g)$ be a smooth Fredholm manifold with oriented Riemannian $q$-structure $(1 \leq q \leq \infty)$. If $M$ has a $\Spin_q$-structure then there are isomorphisms
\begin{equation*}
K^{\infty - j}(M, \F) \cong K_{j+1}(\A(M, g, \F)) \cong K_{j}^c(M), \quad j = 0, 1,
\end{equation*}
where $\F = (M_n , U_n)_{n= k}^\infty$ is {\bf any} augmented Fredholm filtration of $M$.
\end{thm}
\noindent Thus,  the $K$-theory groups of $(M, \F)$ and of the $\cs$-algebra $\A(M, g, \F)$ do not depend on the choice of the Riemannian metric $g$ or the (augmented) Fredholm filtration $\F$.  The dimension shift and the relation with Poincar\'e duality for finite-dimensional spin manifolds then justifies our interpretation of $\A(M, g, \F)$ as an appropriate non-commutative (suspension of the) ``algebra of functions vanishing at infinity'' on $M$.

Finally, it should be noted that, given a Fredholm filtration $\{M_n\}_{n=k}^\infty$ of $M$, we can also naturally associate an {\it inverse} limit algebra, called by Phillips \cite{Phi88} a $\sigma$-$\cs$-algebra, 
\begin{equation*}
C^{\text{inv}}_0(M) = \varprojlim \, C_0(M_n)
\end{equation*}
where the connecting map $C_0(M_{n+1}) \to C_0(M_n)$ is the pullback under the inclusion $M_n \hookrightarrow M_{n+1}$. However, this algebra does {\it not} have the structure of a $\cs$-algebra, in general.  Moreover, if we try to define  the ``topological $K$-theory'' of $M$ as the {\it inverse} limit (using contravariance of topological $K$-theory)
\begin{equation*}
K^j_\text{inv}(M) =  \varprojlim \,  K^j(M_n), \quad j = 0, 1,
\end{equation*}
then we do not get a well-behaved functor. Indeed, as Buhshtaber and Mishchenko have shown, the resulting  $K$-theory sequence of a pair $(M, N)$ is not exact, in general \cite{BuMis2,BuMis1} even for $CW$-complexes. Also, $K$-theory does not behave well with respect to inverse limits since there is a Milnor $\varprojlim \!{}^1$-sequence (Theorem 3.2 \cite{Phi89-2}):
\begin{equation*}
0 \to \varprojlim \!{}^1 K^{j+1}(M_n) \to RK_j(C^{\text{inv}}_0(M)) \to K^j_\text{inv}(M) \to 0
\end{equation*}
where $RK_j$ is the representable $K$-theory for $\sigma$-$\cs$-algebras developed by Phillips \cite{Phi89-2} and Weidner \cite{Weid89}. Hence, there would be no corresponding Serre-Swan duality theorem as in the finite-dimensional category.

The authors would like to thank John Roe, Carolyn Gordon, David Webb, Dana Williams, Gregory Leibon, and the referee for interesting discussions and helpful suggestions.

\section{Clifford $\cs$-algebras and the Thom $*$-Homomorphism}

In this section we assemble the constructions and results for finite-dimensional manifolds that are needed to carry out the direct limit construction of the $\cs$-algebra of an infinite-dimensional Fredholm manifold. All of the manifolds in this section are assumed to be smooth, Hausdorff, paracompact, and finite-dimensional. For a detailed discussion of most of the results in this section, including more proofs, see Section 2 of Trout \cite{Trou03}. 

Let $V$ be a finite-dimensional Euclidean vector space with inner product $\langle \cdot, \cdot \rangle$. The complex Clifford algebra of $V$, denoted $\Cliff(V)$,  is the universal complex $\cs$-algebra (with unit) generated by the elements of $V$ such that
$v^* = v$ and $$v \cdot w + w \cdot v = 2 \langle v, w \rangle 1$$
for all $v, w \in V$. It has a natural $\Z_2$-grading by declaring that all elements of $V$ have odd degree.
The universal property  \cite{LM89,GVF01} of $\Cliff(V)$ is that if $f : V \to A$ is a real linear map of $V$ into a unital complex $\cs$-algebra $A$ such that 
\begin{equation*}
f(v)^2 = \langle v, v \rangle 1_A
\end{equation*}
for all $v \in V$ then there is an induced $\cs$-algebra homomorphism $\tilde{f} : \Cliff(V) \to A$ such that the following diagram commutes:
\begin{equation*}
\xymatrix{ \Cliff(V) \ar[dr]^{\tilde{f}} & \\
V \ar[u]^{C} \ar[r]_{f} & A}
\end{equation*}
where we denote by $C : V \hookrightarrow \Cliff(V)$ the canonical inclusion. However, we will usually identify $v = C(v) \in \Cliff(V)$ for all $v \in V$. 
An important property of these $\Z_2$-graded $\cs$-algebras is their behavior with respect to orthogonal sums:
\begin{equation}\label{eqn:orthogonality}
\Cliff(V \oplus W) \cong \Cliff(V) \gtimes \Cliff(W)
\end{equation}
where $\gtimes$ denotes the $\Z_2$-graded tensor product. $($See the books \cite{LM89,GVF01} for a review of Clifford algebras and Blackadar  \cite{Bla98} for a review of graded $\cs$-algebras.$)$

Let $M_n$ be a {\it finite-dimensional}  smooth Riemannian manifold of dimension $n$ with Riemannian metric $g$. Let $TM_n \to M_n$ denote the tangent bundle of $M_n$.  Let $\Cliff(TM_n) \to M_n$ denote the Clifford bundle \cite{ABS64,BGV86} of $TM_n$, \ie  the bundle of Clifford algebras over $M_n$ whose fiber at
$x \in M_n$ is the complex Clifford algebra $\Cliff(T_xM_n)$ of the Euclidean tangent space $T_xM_n$. It has an induced $\Z_2$-graded $\cs$-algebra bundle structure.

\begin{dfn}\label{def:man-alg}\cite{Kas88} 
Denote by $\CC(M_n)$ the $\cs$-algebra  
\begin{equation*}
\CC(M_n) = C_0(M_n, \Cliff(TM_n))
\end{equation*}
 of continuous sections of $\Cliff(TM_n)$ which vanish at
infinity on $M_n$, with induced $\Z_2$-grading from $\Cliff(TM_n)$. (Kasparov \cite{Kas88} used the notation $C_\tau(M_n)$.)
  \end{dfn}
  
For example, if $M_n = V$ is a finite dimensional Euclidean vector space, then $TM_n \cong V \times V$ and so $\CC(M_n) \cong C_0(V, \Cliff(V))$ as in Definition 2.2  \cite{HKT98}. A priori, this $\cs$-algebra depends on the Riemannian metric $g$ of $M_n$. However, the universal property of Clifford algebras shows that the $\cs$-algebra structure on $\CC(M_n)$ depends only on the manifold $M_n$ and not the chosen metric $g$. Indeed, if $h$ is another Riemannian metric on $M_n$, then $\alpha = \hat{h}^{-1} \circ \hat{g} : TM_n \to TM_n$ is an automorphism of the tangent bundle $TM_n$, where $\hat{g} : TM_n \to T^*M_n$ is the (co)tangent bundle isomorphism induced by any metric $g$. It satisfies
$$h(\alpha(X), X) = g(X, X) \geq 0$$ for any vector field $X$. Thus, $\alpha$ is positive definite with respect to the metric $h$ and so has a positive square root, \ie a bundle automorphism $\beta : TM_n \to TM_n$ such that
$$h(\beta(X), \beta(X)) = h(\alpha(X), X) = g(X, X).$$
If $\Cliff(TM_n, h)$ denotes the Clifford bundle of $M_n$ with respect to the metric $h$ then
$$\beta(X)^2 = g(X, X)1$$
in $\Cliff(TM_n, h)$. By the universal property above (applied to each fiber) $\beta$ extends to an isomorphism
$\tilde{\beta} : \Cliff(TM_n, g) \to \Cliff(TM_n, h)$ of Clifford bundles. (See also Section 9.1 \cite{GVF01}.) By taking sections, there is a canonically induced isomorphism $$\hat{\beta} : \CC(M, g) \to \CC(M, h)$$ of $\Z_2$-graded $\cs$-algebras.

Let $C_0(M_n)$ denote the commutative $\cs$-algebra of continuous complex-valued functions on $M_n$ vanishing at infinity. We always consider $C_0(M_n)$ to be trivially graded. If a $\Z_2$-graded $\cs$-algebra $A$ is equipped with a (fixed)  $*$-homomorphism $\Theta : C_0(M_n) \to Z(M(A))$ that is nondegenerate and has grading degree zero, where $Z(M(A))$ denotes the center of the multiplier algebra of $A$, then we say that $A$ has a $\Z_2$-graded $C_0(M_n)$-algebra structure  \cite{Trou03}. We denote $\Theta(f) a = f \cdot a$ for all $f \in C_0(M_n)$ and $a \in A$. Note that pointwise multiplication
$(fs)(x) = f(x) s(x), \forall x \in M_n,$ where $f \in C_0(M_n)$ and $s \in \CC(M_n)$, determines a nondegenerate $*$-homomorphism $C_0(M_n) \to ZM(\CC(M_n))$  into the center of the multiplier algebra of $\CC(M_n)$ of grading degree zero. Thus, we have the following.

\begin{cor}\label{cor:indep} 
The $\cs$-algebra $\CC(M_n)$ has a canonical $\Z_2$-graded $C_0(M_n)$-algebra structure, and up to $\Z_2$-graded isomorphism, is independent of the Riemannian metric on $M_n$.
\end{cor}

\begin{dfn}\label{def:C*algMn}
Let $\S$ denote the $\cs$-algebra $C_0(\R)$ of continuous complex-valued functions on the real line which vanish at
infinity, with $\Z_2$-grading by even and odd functions. If $A$ is any $\Z_2$-graded $\cs$-algebra then we let $\S A$ be the graded (max) tensor product $\S \gtimes A$.
In particular, let 
\begin{equation*}
\A(M_n) =_{\textrm{def}} \S \CC(M_n) = \S \gtimes C_0(M_n, \Cliff(TM_n))
\end{equation*}
which can be viewed as a non-commutative topological suspension of $M_n$.\footnote{Recall that the suspension of a $\cs$-algebra $A$ is the $\cs$-algebra $SA = C_0(\R) \otimes A$. In particular,  $SC_0(M_n) \cong C_0(\R \times M_n)$ where $\R \times M_n$ is the (reduced) topological suspension of $M_n$.}
\end{dfn}

The following functoriality result will be used when we identify the total space of the normal bundle of an embedding with an open tubular neighborhood. 

\begin{lem}\label{lem:diffeo}\cite{Trou03}
Let $\phi : M_n \to N_n$ be a diffeomorphism of  Riemannian manifolds. There is an induced $\Z_2$-graded $\cs$-algebra isomorphism
 \begin{equation*}
 \phi_* : \A(M_n) \to \A(N_n).
 \end{equation*}
\end{lem} 

\Pf Let $g$ denote the metric on $M_n$ and $h$ denote the metric on $N_n$. If $\phi : M_n \to N_n$ is a diffeomorphism, then, using the pullback metric $\phi^*(h)$, we have that
$$\phi: (M_n, \phi^*(h)) \to (N_n, h)$$
is an {\it isometry} of Riemannian manifolds, which clearly induces a canonical isomorphism
$$\hat{\phi} : \CC(M_n, \phi^*(h)) \to \CC(N_n, h)$$
of $\Z_2$-graded $\cs$-algebras. By the argument above, we have a canonical isomorphism
$$\hat{\beta} : \CC(M_n, g) \to \CC(M_n, \phi^*(h)).$$
Taking the composition and tensoring with the identity of $\S$ gives the required canonical isomorphism
$$\phi_* = id_\S \gtimes  (\hat{\beta} \circ \hat{\phi}) : \A(M_n) = \S \gtimes \CC(M_n) \to \S \gtimes \CC(N_n) = \A(N_n)$$
of $\Z_2$-graded $\cs$-algebras. \hfill \EPf

The following is an easy functoriality property for open inclusions.

\begin{lem}\label{lem:open}\cite{Trou03}
Let $U_n$ be an open subset of the Riemannian manifold $M_n$. The inclusion $i : U_n \into M_n$ induces a short exact sequence
$$\xymatrix{ 0 \ar[r] & \A(U_n) \ar[r]^-{1 \gtimes i_*}  & \A(M_n) \ar[r] &  \A(M_n \setminus U_n) \ar[r] & 0 }$$ 
of $\cs$-algebras. Thus, $\A(U_n) \ideal \A(M_n)$ as a $($two-sided$)$ $\cs$-ideal.
\end{lem}

Let $p : E \to M_n$ be a smooth finite rank Euclidean vector bundle. We will show that there is
a natural ``Thom'' $*$-homomorphism
$$\Psi_p : \A(M_n) \to \A(E),$$  where we consider $E$ as a finite-dimensional manifold with Riemannian structure to be constructed as follows. The main example we have in mind is where $E = \nu M_n$ is the (total space of the) normal bundle of an isometric embedding $M_n \hookrightarrow M_{n+1}$.

Given $p : E \to M_n$, there is a short exact sequence \cite{AMR88,BGV86} of real vector bundles 
$$\xymatrix{ 0 \ar[r] & VE \ar[r] & TE \ar[r]^-{T^*p} & p^*TM_n \ar[r] & 0}$$  
where the {\it vertical subbundle} $VE = \ker(T^*p)$ is isomorphic to $p^*E$. This
sequence does {\it not} have a canonical splitting, in general, but choosing a compatible connection $\nabla$ on $E$ determines an associated vector bundle
splitting. Recall that a connection 
$\nabla : C^\infty(M_n, E) \to C^\infty(M_n, T^*M_n \otimes E)$  
on $E$ is {\it compatible} \cite{BGV86,LM89} with the bundle metric $( \cdot, \cdot )$ on $E$ if 
$$d(s_1,s_2) = (\nabla s_1, s_2) + (s_1, \nabla s_2)$$  for all smooth sections $s_1, s_2 \in C^\infty(M_n,E)$.
If $p: E \to M_n$ is equipped with a compatible connection $\nabla$, then we call $E$ an {\bf affine} Euclidean bundle. 

Let $\nabla^{*} :C^\infty(E, p^*E) \to C^\infty(E, T^*E \otimes p^*E)$ denote the pullback of $\nabla$ on the bundle $p^*E \to E$, which is defined by the formula:
$$\nabla^{*}(f p^*s) = df \otimes p^*s + f p^*(\nabla s)$$ for $f \in C^\infty(M_n)$ and $s \in C^\infty(M_n,E)$. The {\it tautological section} $\tau \in C^\infty(E, p^*E)$ is the smooth section of $p^*E \to E$ defined by the formula $\tau(e) = (e, e)$ for all $e \in E$. The derivative of $\tau$ will be denoted by
$$\omega = \nabla^{*} \tau \in C^\infty(E, T^*E \otimes p^*E) = \Omega^1(E, p^*E) \cong \Omega^1(E, VE)$$ 
which is the connection $1$-form of $\nabla$ (see Definition 1.10 \cite{BGV86}). The kernel $HE = \ker(\omega) \cong p^*TM_n$ of the connection $1$-form $\omega$ is the horizontal subbundle of $TE$ which provides a splitting
$$TE = VE \oplus HE \cong p^*E \oplus p^*TM_n.$$  Now give $TE$ the direct sum of the pullback metrics on $p^*E$ and $p^*TM_n$. This gives $E$ the structure of a Riemannian manifold and makes the splitting of $TE$ orthogonal. 

\begin{lem}\label{lem:splitting}\cite{Trou03}
Let $p: E \to M_n$ be a finite rank affine Euclidean bundle on the Riemannian manifold $M_n$. There is an induced orthogonal splitting of the exact sequence
$$0 \to p^*E \to TE \to p^*TM_n \to 0$$ and so  there is a canonical isomorphism of Euclidean vector bundles
$$TE \cong p^*E \oplus p^*TM_n$$  where $p^*E$ and $p^*TM_n$ have the pullback metrics. Thus ,the manifold $E$ has a canonical Riemannian metric.
\end{lem}


Hence, given a compatible connection $\nabla$ on the Euclidean bundle $E$, we can define the $\cs$-algebra $\CC(E)$ as above using the induced
Riemannian structure on the {\it manifold} $E$. However, we also have the $\cs$-algebra $C_0(E, \Cliff(p^*E))$ associated to the pullback bundle
$p^*E \to E$.\footnote{Note: Although $C_0(E, \Cliff(p^*E)) \cong p^*C_0(M_n, \Cliff(E))$,
we will not need this isomorphism.} Both $\CC(E)$ and $C_0(E, \Cliff(p^*E))$ have natural
$C_0(E)$-algebra structures. However, the bundle map $p: E \to M_n$ induces a pullback
$*$-homomorphism \cite{RW85} $$p^* : C_0(M_n) \to C_b(E) = M(C_0(E)) : f \mapsto p^*(f) = f \circ p$$  which induces a (graded) $C_0(M_n)$-algebra structure on any
(graded) $C_0(E)$-algebra.  

\begin{dfn} \cite{Trou03} Let $A$ and $B$ be $\Z_2$-graded $C_0(M_n)$-algebras. The balanced tensor product over $M_n$, denoted $A \gtimes_{M_n} B$, is the quotient of the maximal graded tensor product $A \gtimes B$ \cite{Bla98} by the ideal $J$ generated by 
\begin{equation*}
\{(f \cdot a) \gtimes b - a \gtimes (f \cdot b) : a \in A, b \in B, f \in C_0(M_n)\}.
\end{equation*}
For example $C_0(M_n) \gtimes_{M_n} A \cong A$ via the map induced by  $f \gtimes a \mapsto f \cdot a$.
\end{dfn}

The following is an important result that relates these two $\cs$-algebras to the $\cs$-algebra $\CC(M_n)$ of the base manifold $M_n$.

\begin{thm}\label{thm:decomp} 
Let $p: E \to M_n$ be a finite rank affine Euclidean bundle on the Riemannian manifold $M_n$.
There is a natural isomorphism of graded $\cs$-algebras
$$\CC(E) \cong C_0(E, \Cliff(p^*E)) \gtimes_{M_n} \CC(M_n).$$ 
\end{thm}

\Pf By the previous lemma, there is an induced orthogonal splitting 
$$TE = p^*E \oplus p^*TM_n.$$   Thus,  we have an induced isomorphism of $\Z_2$-graded Clifford algebra bundles
\begin{eqnarray}\label{eqn:cliffbdl}
\Cliff(TE) \cong \Cliff(p^*E \oplus p^*TM_n) = \Cliff(p^*E) \gtimes p^*\Cliff(TM_n).
\end{eqnarray} Therefore, by taking sections, we have canonical balanced tensor product isomorphisms (see Proposition A.7 \cite{Trou03})
\begin{align*}
\CC(E) & =_{\text{def}}   C_0(E, \Cliff(TE)) \cong C_0(E, \Cliff(p^*E) \gtimes p^*\Cliff(TM_n)) \\
       & \cong_{\phantom{\text{def}}}  C_0(E, \Cliff(p^*E)) \gtimes_{_E} C_0(E, \Cliff(p^*TM_n)).
\end{align*} But, we have, using pullbacks along $p : E \to M_n$, that there are canonical pullback isomorphisms 
(see Proposition A.9  \cite{Trou03})
$$C_0(E, \Cliff(p^*TM_n)) \cong p^*C_0(M, \Cliff(TM_n)) = p^*\CC(M_n) =_{\text{def}} C_0(E) \gtimes_{M_n} \CC(M_n).$$ 
Hence, it follows that
\begin{align*}
\CC(E) & \cong  C_0(E, \Cliff(p^*E)) \gtimes_{_E} C_0(E, \Cliff(p^*TM_n)) \\
       & \cong  C_0(E, \Cliff(p^*E)) \gtimes_{_E} C_0(E) \gtimes_{M_n} \CC(M_n) \\
       &  \cong  C_0(E, \Cliff(p^*E)) \gtimes_{M_n} \CC(M_n)
\end{align*} using the canonical isomorphism $A \gtimes_{_E} C_0(E) \cong A$ for graded $C_0(E)$-algebras. \EPf 

We now wish to define a certain ``Thom operator'' for the ``vertical'' algebra $C_0(E, \Cliff(p^*E))$. Associate to the Euclidean bundle $E$ an unbounded section
\begin{equation*} C_E : E  \to  \Cliff(p^*E) : e  \mapsto C_{p(e)}(e) \end{equation*}   where $C_{p(e)}$ is the Clifford operator on the Euclidean
space $E_{p(e)}$ from Definition 2.4 of \cite{HKT98}. It is given globally by the composition
$$\xymatrix{ E \ar[r]^-{\tau} \ar@/_1.5pc/[rr]_-{C_E}  & p^*E \ar[r]^-{C} & \Cliff(p^*E)}$$   where $\tau \in C^\infty(E, p^*E)$ is the tautological section
(see above) and
$C : p^*E \into \Cliff(p^*E)$ is the canonical inclusion $C(e_1, e_2) = C_{p(e_1)}(e_2)$.   The following is then easy to prove.

\begin{thm}\label{thm:cliffop} 
Let $E$ be a finite rank Euclidean bundle on $M_n$. Multiplication by the section $C_E : E \to \Cliff(p^*E)$  determines a
degree one, essentially self-adjoint, unbounded multiplier $($see Definition A.1 \cite{Trou03}$)$ of the $\cs$-algebra $C_0(E, \Cliff(p^*E))$ with domain $C_c(E, \Cliff(p^*E))$.
\end{thm}

We will call $C_E$ the {\bf Thom operator} of $E \to M_n$. Thus, we have a functional calculus homomorphism
$$\S  \to M(C_0(E, \Cliff(p^*E))) : f \to f(C_E)$$  from $\S$ to the multiplier algebra of $C_0(E, \Cliff(p^*E))$. Note that $f(C_E)$ goes to
zero in the ``fiber'' directions on $E$ (since $p(e)$ is constant), but is only bounded in the ``manifold'' directions on $E$.  Indeed, for
the generators $f(x) = \exp(-x^2)$ and $g(x) = x \exp(-x^2)$ of $\S$, we have that $f(C_E)$ and $g(C_E)$ are, respectively, multiplication by the
following functions on $E$:
$$f(C_E)(e) = \exp(-\|e\|^2) \text{ and } g(C_E)(e) = e \cdot \exp(-\|e\|^2), \quad \forall e \in E.$$ 

\begin{dfn} Let $X$ denote the degree one, essentially self-adjoint, unbounded multiplier of $\S$, with domain the compactly supported functions, given by multiplication by $x$, \ie $Xf(x) = xf(x)$ for all $f \in C_c(\R)$ and $x \in \R$. \end{dfn}

By Lemma A.3 \cite{Trou03}, the operator $X \gtimes 1 + 1 \gtimes C_E$ determines a degree one, essentially self-adjoint, unbounded multiplier of 
the tensor product
$$\S \gtimes C_0(E, \Cliff(p^*E)) = \S C_0(E, \Cliff(p^*E))$$ with domain $C_c(\R) \hat{\odot} \, C_c(E, \Cliff(p^*E))$. We obtain a functional calculus homomorphism 
$$\beta_E : \S \to M(\S C_0(E, \Cliff(p^*E))) : f \mapsto f(X \gtimes 1 + 1 \gtimes C_E)$$
from $\S$ into the multiplier algebra of $\S C_0(E, \Cliff(p^*E))$.  Now we can define our ``Thom $*$-homomorphism'' for a
finite rank affine Euclidean bundle. This will provide part of the connecting map in Section 4 when we define the direct limit $\cs$-algebra for an infinite-dimensional Riemannian Fredholm manifold.

\begin{thm}\label{thm:Thomhom} 
Let $p: E \to M_n$ be a finite rank affine Euclidean bundle on the Riemannian manifold $M_n$. With respect to the isomorphism
$$\A(E) \cong \S \gtimes C_0(E, \Cliff(p^*E)) \gtimes_{M_n} \CC(M_n)$$ from Theorem \ref{thm:decomp}, there is a graded $*$-homomorphism 
$$\Psi_p  =  \beta_E \gtimes_{M_n} \id_{M_n} : \A(M_n) \to \A(E)$$    which on elementary tensors $f \gtimes s \in \S \gtimes \CC(M_n)=\A(M_n)$ is given by
$$f \gtimes s \mapsto f(X \gtimes 1 + 1 \gtimes C_E) \gtimes_{M_n} s.$$
\end{thm}

\Pf From the discussion above, we have that $\beta_E \gtimes_{M_n} \id_{M_n}$ is the composition
$$\xymatrix{
\A(M_n) \ar[r]^-{\beta_E \gtimes \id} & M(\S C_0(E, \Cliff(p^*E))) \gtimes \CC(M_n) \to  M(\S C_0(E, \Cliff(p^*E))) \gtimes_{M_n} \CC(M_n)
}$$
Checking on the generator $f(x) = \exp(-x^2)$ of $\S$, we compute that
$$f(X \gtimes 1 + 1 \gtimes C_E) \gtimes s = \exp(-x^2) \gtimes \exp(-\|e\|^2) \gtimes_{M_n} s \in \A(E)$$ Similarly for $g(x) = x \exp(-x^2)$,
we find that
\begin{eqnarray*} g(X \gtimes 1 + 1 \gtimes C_E) \gtimes_{M_n} s & =  &x\exp(-x^2) \gtimes \exp(-\|e\|^2) \gtimes_{M_n} s \\
 & + & \exp(-x^2) \gtimes e \cdot \exp(-\|e\|^2)\gtimes_{M_n} s\in \A(E).
\end{eqnarray*}  It follows that the range of $\Psi_p = \beta_E \gtimes_{M_n} \id_{M_n}$ is in $\A(E)$ as desired.
\EPf

Since the space of compatible connections $\nabla$ on $E \to M_n$ is convex, we have the following result.

\begin{prop}\label{prop:htpy} 
Let $p : E \to M_n$ be a smooth finite rank affine Euclidean bundle  on the Riemannian manifold $M_n$. The homotopy class of the  $*$-homomorphism $\Psi_p : \A(M_n) \to
\A(E)$ is independent of the choice of compatible connection on $E$. 
\end{prop}

\begin{prop}\label{prop:triv} 
If $p : E = M_n \times V \to M_n$ is a trivial finite rank affine Euclidean bundle $($with trivial connection $\nabla_0 = d$ $)$ then we
have a $\Z_2$-graded isomorphism
$$\CC(E) \cong \CC(V) \gtimes \CC(M_n)$$ such that the Thom map has the form
$$\Psi_p \cong \beta_V \gtimes \id_{\CC(M_n)} : \A(M_n) = \S \gtimes \CC(M_n) \to \A(V) \gtimes \CC(M_n) \cong \A(E)$$   where $\beta_V : \S \to \A(V) : f \mapsto
f(X \gtimes 1 + 1 \gtimes C_V)$ is the Thom map for $V \to \{0\}$.
\end{prop}

\Pf The trivial connection $\nabla_0 = d$ gives the manifold $E = M_n \times V$ the Riemannian metric induced by the isomorphism 
$$TE = TM_n \times TV \to M_n \times V =  E.$$ The pullback vector bundle $p^*E \to E$ has the form
$$p^*E = (M_n \times V) \times V \to M_n \times V = E$$ and so the Clifford bundle $\Cliff(p^*E) = (M_n \times V) \times \Cliff(V)$, which gives:
$$C_0(E, \Cliff(p^*E)) = C_0(M_n \times V, (M_n \times V) \times \Cliff(V)) \cong C_0(V, \Cliff(V)) \gtimes C_0(M_n).$$ 
By Theorem \ref{thm:decomp}, it
follows that
\begin{align*}
\CC(E) & \cong  C_0(E, \Cliff(p^*E)) \gtimes_{M_n} \CC(M_n) \\
       & \cong  C_0(V, \Cliff(V)) \gtimes C_0(M_n) \gtimes_{M_n} \CC(M_n) \\
       & \cong  \CC(V) \gtimes \CC(M_n).
\end{align*} 
where we used the isomorphism $C_0(M_n) \gtimes_{M_n} \CC(M_n) \cong \CC(M_n)$. The result now easily follows. \EPf

For example, if $p: E_b \to E_a$ is the orthogonal projection of a finite dimensional Euclidean vector space $E_b$ onto a linear subspace $E_a$ then
$\Psi_p = \beta_{ba}$ is the  ``Bott homomorphism'' from Definition 3.1 of Higson-Kasparov-Trout \cite{HKT98}. 

\section{Fredholm Manifolds and Filtrations}

\newcommand{\bdd}{\mathcal{L}}
\newcommand{\gd}{\delta}
\newcommand{\gm}{\gamma}
\newcommand{\lan}{\langle}
\newcommand{\ra}{\rightarrow}
\newcommand{\ran}{\rangle}


Fredholm manifolds are a particular case of {\it Hilbert manifolds},
\ie manifolds modeled on a separable infinite-dimensional real Hilbert space.
Most of the standard constructions from the differential geometry of
finite-dimensional manifolds carry on in the infinite dimensional
situation (as reference see Lang's book \cite{Lang}). 
%
%
All the Hilbert manifolds that we 
consider in this paper are assumed to be
connected, separable, paracompact, Hausdorff, and infinitely smooth.

Let $\E$ be a separable infinite-dimensional Euclidean space, \ie a real
Hilbert space of countably infinite dimension. We will use the following 
notation:
$\bdd(\E)$ denotes the real $\cs$-algebra of bounded linear operators on $\E$;
$\F=\F(\E)$ denotes the finite rank operators;
$\K=\K(\E)$ denotes the closed ideal of compact operators;
$\Phi=\Phi(\E)$ denotes the Fredholm operators;
and $GL(\E)$ denotes  the Banach-Lie group of units of $\bdd(\E)$,
with identity $I$.

\begin{dfn}\label{D:perturbation}
A \emph{perturbation class $P$ of $\E$} is a subspace $P=P(\E)$ of
$\bdd(\E)$ such that: (1) $\F(\E)\subseteq P(\E)$, 
(2) $P(\E)$ is an ideal in $\bdd(\E)$, and
(3) $\Phi(\E) + P(\E) = \Phi(\E)$.
\end{dfn}

\noindent
As examples of perturbation classes we have: the finite rank operators $\F(\E)$, the compact operators $\K(\E)$, or indeed any proper two-sided ideal included in $\K$. For $1 \leq q < \infty$, let $P_q$ be the perturbation class defined as the closure of $\F(\E)$ under the norm 
\begin{equation}
\|\, T\,\|_q =\left( \text{Trace}(T^*T)^{q/2} \right)^{1/q}.
\end{equation}
If $q=1$ one obtains the trace-class operators, and if $q=2$ the
Hilbert-Schmidt operators. If $q = \infty$, then we set $P_\infty = \K(\E)$ with norm $\|\, T \|_\infty = \|\, T \|$. 

Given a perturbation class $P$ of $\E$, we let
$$
\begin{aligned}
GL_P(\E) &=  GL(\E) \cap (I+P(\E)) \\
         &= \{ T=I+K \;|\; T\in GL(\E), K\in P(\E) \}.
\end{aligned}
$$
For $1 \leq q < \infty$, we  abbreviate $GL_{P_q}(\E) = GL_q(\E)$.  For $p = \infty$, we abbreviate $GL_{\K(\E)}(\E) = GL_\K(\E)$. We topologize $GL_q(\E)$ by requiring that the map $GL_q(\E) \to \O \subset P_q$,  $K \mapsto I+K,$ be a  homeomorphism, where $\O$ is the set of all $K$ with $I+K$ invertible \cite{Pal65b}. In general, $GL_P(\E)$ is a normal subgroup of $GL(\E)$, but, where $GL(\E)$ is contractible (by Kuiper's theorem \cite{Kui65}), $GL_P(\E)$ may not be contractible. For example, by Theorem B of Palais \cite{Pal65b}, we have that
$$\pi_0(GL_q(\E))=\Z/2\Z$$
for all $1 \leq q \leq \infty$. However, $GL_P(\E)$ is not a {\it closed} subgroup of $GL(\E)$ unless $P = P_\infty = \K(\E)$.

\begin{dfn}\label{def:fredmfld}
Let $M$ be a Hilbert manifold modeled on $\E$. A \emph{Fredholm structure on $M$}
is an integrable reduction of the principal $GL(\E)$-bundle of $M$ to 
$GL_{\K}(\E)$. Equivalently, it is a maximal atlas of $M$ such that the 
differential of the change of coordinates maps is an element of $GL_{\K}(\E)$ 
at every point.
A \emph{Fredholm manifold} is a Hilbert manifold with a 
specified Fredholm structure.
\end{dfn}

Since there is a natural inclusion  $GL_P(\E) \hookrightarrow GL_\K(\E)$ induced by the inclusion $P(\E) \hookrightarrow \K(\E)$, if a Hilbert manifold $M$ is equipped with a reduction of it's structure group from $GL(\E)$ to $GL_P(\E)$ then we can give $M$ a canonical Fredholm structure in the sense of the previous definition. We will make use of this fact when discussing spin structures for Fredholm manifolds in Section 5.

\noindent
{\bf Note.} A $C^{\infty}$-map $f : M \ra N$ between Hilbert manifolds is called 
a {\it Fredholm map} if, for every $x\in M$, $Df(x):T_xM \ra T_{f(x)}N$
is a Fredholm operator. Fredholm manifolds are exactly the manifolds on which 
Fredholm maps can be constructed. Results of Elworthy and Tromba 
\cite{ElwTr}
show that for 
a Fredholm manifold $M$ there is an index zero (even {\it bounded} and {\it proper}) 
Fredholm map $f : M \ra \E$.

The following decomposition theorem is crucial in the study of Fredholm
manifolds (\cite[Thm 2.2]{Mkhr1}):
\begin{thm}\label{Thm:FredFilt}
Let $M$ be a Fredholm manifold. There exists a sequence $\{M_n\}_{n=k}^{\infty}$
of finite dimensional closed submanifolds such that:
\begin{itemize}
\item[(i)]   $\text{dim } M_n = n$; $M_n \subset M_{n+1}$;
\item[(ii)]  the inclusions $M_n \hookrightarrow M_{n+1}$ and
$M_n \hookrightarrow M$ have trivial normal bundles;
\item[(iii)] $M_{\infty} = \cup_{n\geq k} M_n$ is dense in $M$; and
\item[(iv)]  the natural inclusion map $M_{\infty} \hookrightarrow M$ is a 
homotopy equivalence,  if $M_{\infty}$ is given the {\rm direct limit topology}.
\end{itemize}
\end{thm}

\noindent
A sequence $\{M_n\}_{n=k}^{\infty}$ as in the theorem above is called a
\emph{Fredholm filtration of $M$}.
%

We will now give some examples (and a non-example) of Fredholm manifolds and filtrations. 

\begin{examples}
(i)~The Euclidean space $M=\E$ has an obvious Fredholm structure, determined
by a single chart $I:\E\ra\E$. It is the only possible structure.
Let $\{e_n\}_{n=1}^{\infty}$ be an orthonormal basis of $\E$, and 
$E_n$ be the linear span of $\{e_1, e_2, \dots, e_n\}$. The sequence
$\{ E_n\}_n$ is known as a {\it flag} of $\E$, and it forms 
a Fredholm filtration. 
\newline (ii)~The unit sphere of $\E$, $S_{\E}=\{ x\in \E \,|\, \|x\|=1\}$, gets
by restriction from $\E$ a Fredholm structure. As a Fredholm filtration we have
$$S^1 \subset S^2 \subset \dots \subset S^n \subset \dots \subset S_{\E}.$$
(iii) The following is a {\it non-example}.
The sequence of real projective spaces
$$
\R P^1 \subset \R P^2 \subset \dots \subset \R P^n \subset \dots \subset \R P_{\E}
$$
is not a Fredholm filtration of the infinite dimensional real projective space $\R P_\E$ of $\E$, for any choice of Fredholm structure, because the inclusions $\R P^n \subset \R P^{n+1}$ do not have trivial normal bundles.
\end{examples}

To get an idea how Fredholm filtrations are constructed in general, we briefly outline the procedure as follows. Let $M$ be a Fredholm manifold modeled on $\E$. Let $\{E_n\}_n$ be a flag for $\E$ as in Example 3.4 (i) above.
Choose an index zero Fredholm map $f : M \ra \E$ which is transversal to the $E_n$'s, and define 
$M_n = f^{-1} (E_n)$. Each $M_n$ (when nonempty) is a finite-dimensional submanifold of $M$ of dimension $n$ and $M_n \subset M_{n+1}$. The normal bundle $\nu M_n \to M_n$ of the inclusion $M_n \subset M_{n+1}$ is the pullback $\nu M_n = f^*(\nu E_n)$ of the (trivial) normal bundle $\nu E_n = E_n^\perp \cap E_{n+1}$ and, hence, is trivial. The sequence $\{M_n\}_{n=k}^\infty$, where $M_k \neq \emptyset$ is the first nonempty submanifold, forms a Fredholm filtration of $M$. Note that since there is always a {\it bounded, proper} index zero Fredholm map $f : M \to \E$, the $M_n$'s can be chosen to be {\it compact}. See the Addendum to Theorem 2C in Eells and Elworthy \cite{EeElw2}.

One can actually say more about the Fredholm filtrations of a Fredholm manifold,
but we need to recall first some facts about the differential geometry of infinite
dimensional manifolds.

\begin{dfn}\label{D:tubularN}
Let $N$ be a submanifold of $M$. A {\it tubular neighborhood
of $N$ in $M$} consist of the following data:  a vector bundle
$\pi: B \rightarrow N$ over $N$, an open neighborhood $V$ of the 
zero section $\zeta(N)$ in $B$, an open set $U$ in $M$ containing $N$,
and a diffeomorphism 
$f:V \rightarrow U$ which commutes with the zero section $\zeta:N\ra V$:
\begin{equation*}
\xymatrix{ V \ar[d]^{\pi|_V} \ar[dr]^{f} & \\
N \ar@/^/[u]^{\zeta} \ar@{^{(}->}[r]_{i} & U}
\end{equation*}
$U$ is called the {\it tube} of the tubular neighborhood. The tubular neighborhood is called {\it total} if $V = B$ the total space of the bundle.
\end{dfn}

\noindent Using the notion of 
spray \cite[IV.3]{Lang}, its associated exponential map, and restriction
to the normal bundle of the inclusion $i: N \rightarrow M$, one can
prove the existence and uniqueness of tubular neighborhoods, if $M$ is 
a Hilbert manifold (\cite{Lang}, Theorems IV.5.1 and IV.6.2). On a Riemannian manifold one can always choose tubular neighborhood to be total. 

\begin{dfn}
A {\it Riemannian manifold} is a pair $(M, g)$, where $M$ is a
Hilbert manifold, and $g$ is a metric on $M$, \ie $g_x$ is a 
(smoothly varying) positive-definite non-singular symmetric bilinear form on $T_xM$, for
every $x\in M$.
\end{dfn}

\noindent
According with \cite[Cor.II.3.8]{Lang}, every paracompact 
$C^{\infty}$-manifold modeled on a separable Hilbert
space admits partitions of unity of class $C^{\infty}$.
It follows that Hilbert manifolds admit Riemannian metrics:

\begin{prop}
\cite[Prop.VII.1.1]{Lang}
Let $M$ be a manifold admitting partitions of unity, and let 
$\pi: B \rightarrow M$ be a vector bundle whose fibers are Hilbertable 
vector spaces. Then $\pi$ admits a Riemannian metric.
\end{prop}

Granted all of this, the next statement is a combination of 
\cite[Thm 2.3]{Mkhr1} and remarks from \cite{Mkhr2} and \cite{EeElw2}.

\begin{thm}\label{T:aFf} 
Let $M$ be a Fredholm manifold with Riemannian metric $g$ compatible with the topology of $M$. There exists a Fredholm filtration
$\{M_n\}_{n=k}^{\infty}$ of $M$ for which geodesically defined exponential 
neighborhoods $Z_n$ of $M_n$ in $M$ can be constructed satisfying:
$$
Z_n\subset Z_{n+1} \text{  and  } \cup_{n\geq k} Z_n = M.
$$
Moreover $U_n = Z_n \cap M_{n+1}$ is a tubular neighborhood of 
$M_n$ in $M_{n+1}$, for each $n\geq k$.
\end{thm}

\begin{dfn}\label{D:aFf} 
We call a Fredholm filtration $\{M_n\}_{n=k}^{\infty}$ together with
a collection $\{U_n\}_{n=k}^\infty$, where $U_n$ is a {\it total} tubular neighborhood of $M_n \subset M_{n+1}$,  
an {\it augmented Fredholm filtration} and we shall denote this by 
$\F = (M_n, U_n)_{n=k}^\infty$. Note that we assume that each $U_n$ is equipped with a fixed diffeomorphism $\phi_n : \nu M_n \to U_n$.
\end{dfn}

\noindent
%


Fredholm manifolds often arise as spaces of paths and we end this section with
one more example. In Section 5, Example 5.10, we will discuss examples of Fredholm manifolds arising from loop groups $\Omega G$ of certain compact Lie groups $G$ (and their associated spin structures.)

\begin{example}
~See \cite{EeElw1}.
Let $X$ be a complete finite-dimensional Riemannian manifold, and $a\in X$.
Let $M=P_a(X)$ be the space of paths $\gm : [0,1] \ra X$, with
$\gm(0)=a$ and $\gm$ absolutely continuous with square integrable
derivative. Then $M$ is  a separable smooth Hilbert manifold.
Moreover a complete Riemannian structure on $M$ is given by
$$
g_\gamma(u, v) = \lan u,v \ran_\gm = \int_0^1 \lan D_\gm u, D_\gm v \ran_\gm,
$$
for $u, v \in T_\gm M$, where $D_\gm$ denotes the covariant derivative
along $\gm$. There is  natural diffeomorphism
$$
\gd : P_a(X) \ra P_0(T_aX), 
 \gd(\gm)(t) = \int_0^t \tau_0^s \gm'(s) \,ds,
$$
where $\tau_0^s$ denotes parallel transport along $\gm$ from
$T_{\gm(s)}X$ to $T_aX$. This map $\gd$, called E.~Cartan's development map,
gives a diffeomorphism of $M = P_a(X)$ with the Hilbert space $P_0(T_aX)$ and, hence, a  unique Fredholm structure on the contractible space $M$.
\end{example}

\section{The $\cs$-algebra of a Fredholm Manifold}
 
Let $M$ be a smooth, separable, connected, paracompact Hilbert manifold modeled on the separable,
infinite-dimensional Euclidean space $\E$. We assume that $M$ is equipped with a
Riemannian Fredholm structure, \ie  a reduction of the structure group of $M$
from $GL(\E)$ to $GL_\K(\E)$ and a Riemannian metric $g$ that is compatible with the topology of $M$. This is equivalent to a reduction of the structure group from $GL(\E)$ to $\O_\K(\E) = GL_\K(\E) \cap O(\E)$. (See section 5.2).
 
Let $\F = (M_n, U_n)_{n=k}^\infty$ be an augmented Fredholm filtration of $M$ by closed
$n$-dimensional submanifolds $M_n$ with total tubular neighborhoods $M_n \subset U_n \subset M_{n+1}$, as in Definition \ref{D:aFf}.  Let $p_n : \nu M_n \to M_n$ denote the normal bundle of the embedding $j_n : M_n \hookrightarrow M_{n+1}$. That is, we have a short exact sequence
\begin{equation*}
0 \ra TM_n \ra TM_{n+1}|_{M_n} \ra \nu M_n \ra 0
\end{equation*}
of finite rank vector bundles. 

These geometric considerations lead us to the following topological diagram of bundles and spaces:
\begin{equation}\label{eqn:top.dirlim}
\xymatrix{ 
\nu M_n \ar[d]_{p_n}^{\text{normal}} \ar[rr]^{\phi_n}_{\text{diffeo}} & & 
        U_n \ar@{^{(}->}^{k_n}_{\text{open}}[rr] & & M_{n+1} \\
M_n
}
\end{equation}
where the tubular neighborhood $U_n$ is identified with the total space of the normal bundle $\nu M_n$ via a fixed diffeomorphism $\phi_n : \nu M_n \to U_n$ and $k_n : U_n \hookrightarrow M_{n+1}$ denotes the (open) inclusion.

For each $n$, let $M_n$ have the induced Riemannian metric $g_n = i_n^*(g)$ where $i_n : M_n \hookrightarrow
M$ denotes the inclusion.
Thus, for each $n \geq k$, we have the associated $\cs$-algebra
\begin{equation*}
\A(M_n) = \S  \CC(M_n) = \S \gtimes C_0(M_n, \Cliff(TM_n))
\end{equation*}
as in Definition \ref{def:C*algMn}. Recall that $\S$ denotes the $\cs$-algebra $C_0(\R)$ graded by even and odd functions.

The restricted bundle $TM_{n+1}|_{M_n}$ is the pullback bundle $j_n^*(TM_{n+1})$ under the inclusion $j_n : M_n \hookrightarrow M_{n+1}$. Thus, there is an induced pullback metric $j_n^*(g)$ and pullback connection $j_n^*(\nabla^{n+1})$ on $TM_{n+1}|_{M_n}$, where $\nabla^{n+1}$ is the Levi-Civita connection of $M_{n+1}$ \cite{BGV86}. Using this pullback metric we have an orthogonal splitting
\begin{equation*}
TM_{n+1}|_{M_n} \cong TM_n \oplus \nu M_n
\end{equation*}
of vector bundles on $M_n$. Give $\nu M_n$ the induced bundle metric and projected connection $\nabla^{\nu M_n}$. Thus, $p_n :  \nu M_n \to M_n$ has a canonical structure as an affine Euclidean bundle. By Theorem \ref{thm:Thomhom}, there is an induced $\cs$-algebra homomorphism
\begin{equation*}
\Psi_{p_n} : \A(M_n) \to \A(\nu M_n)
\end{equation*}
where $\nu M_n$ is given the Riemannian metric from Lemma \ref{lem:splitting}.

Give the open set $U_n \subset M_{n+1}$ the induced Riemannian metric $k_n^*(g_{n+1})$ from $M_{n+1}$. By Lemma \ref{lem:open} we have an inclusion of $\cs$-algebras
\begin{equation*}
\xymatrix{(k_n)_* : \A(U_n) \ar@{^{(}->}[r] &   \A(M_{n+1})}
\end{equation*}
induced by the inclusion $k_n : U_n \hookrightarrow M_{n+1}$.
Finally, we have by Lemma \ref{lem:diffeo}, a canonical $\cs$-algebra isomorphism
\begin{equation*}
\xymatrix{(\phi_n)_* : \A(\nu M_n) \ar[r]^-{\cong} & \A(U_n)}
\end{equation*}
induced by the diffeomorphism $\phi_n : \nu M_n \to U_n$ of the tubular neighborhood $U_n$ with the total space $\nu M_n$ of the normal bundle.

Thus, we have the following diagram of $\cs$-algebras and $*$-homomorphisms, which can be considered as the non-commutative version of diagram (\ref{eqn:top.dirlim}) above:
\begin{equation}\label{eqn:alg.dirlim}
\xymatrix{ 
{\A}(\nu M_n) \ar[r]^{(\phi_{n})_*}_{{\cong}} & {\A}(U_n) \ar@{^{(}->}[r]^{(k_n)_*} & {\A}(M_{n+1}) \\
{\A}(M_n) \ar[u]^{\Psi_{_{p_n}}}_{\text{Thom}} \ar@{-->}[urr]_{\alpha_n}
}
\end{equation}
The dotted arrow, which is by definition the composition of the other three, gives the connecting map $\alpha_n : \A(M_n) \ra \A(M_{n+1})$ in the definition of our $\cs$-algebra $\A(M, g, \F)$.

\begin{dfn}\label{def:C*algFM}
Let $M$ be a smooth Fredholm manifold\footnote{Recall that we assume $M$ to be connected, separable, paracompact, and Hausdorff.}, modeled on the separable infinite-dimensional Euclidean space $\E$, equipped with a Riemannian metric~$g$ compatible with the topology of $M$, and an augmented Fredholm filtration $\F = (M_n, U_n)_{n=k}^\infty$.
The {\it $\cs$-algebra of the triple $(M,g,\F)$}
is the direct limit $\cs$-algebra 
\begin{equation}\label{eqn:A(M)}
\A(M,g,\F ) = \varinjlim \, \A(M_n)
\end{equation}
where the direct limit is taken over the directed system $\{\A(M_n),\alpha_n\}_{n=k}^\infty$ and the connecting maps $\alpha_n$ are given by diagram (\ref{eqn:alg.dirlim}).
\end{dfn}
\noindent It easily follows that $\A(M, g, \F)$ has the structure of a  $\Z_2$-graded, separable, nuclear  $\cs$-algebra. One can also show (using Lemma \ref{lem:diffeo} and the construction in Lemma \ref{lem:splitting}) that $\A(M, g, \F)$ does not depend, up to isomorphism of $\Z_2$-graded $\cs$-algebras, on the choice of the Riemannian metric $g$ of $M$. Indeed, we have:

\begin{lem} Let $M$ be a smooth Fredholm manifold with augmented Fredholm filtration $\F = (M_n, U_n)_{n=k}^\infty$. If $g$ and $h$ are Riemannian metrics on $M$ compatible with the topology, there is a canonical map $$\Phi : \A(M, g, \F) \to \A(M, h, \F)$$ which is an isomorphism of $\Z_2$-graded $\cs$-algebras.
\end{lem}

\Pf The identity map $\id_M : (M, g) \to (M, h)$ is a diffeomorphism of Riemannian Fredholm manifolds and  induces for each  $n \geq k$ a commuting diagram
\begin{equation*}
\xymatrix{
\A(M_n, g_n) \ar[r] \ar[d]^{\cong} & \A(\nu M_n, g_n') \ar[r] \ar[d]\ar[d]^{\cong} & \A(U_n, k_n^*(g_{n+1})) \ar[r] \ar[d]\ar[d]^{\cong} & \A(M_{n+1},g_{n+1}) \ar[d]\ar[d]^{\cong}\\
\A(M_n, h_n) \ar[r]           & \A(\nu M_n, h_n') \ar[r]           & \A(U_n, k_n^*(h_{n+1})) \ar[r]          & \A(M_{n+1},h_{n+1})    
}
\end{equation*}        
where $g_n'$ and $h_n'$ are the Riemannian metrics induced on the total space $\nu M_n$ by Lemma \ref{lem:splitting} and the vertical maps are the $\Z_2$-graded $\cs$-algebra isomorphisms induced by $\id_{M_n} : (M_n, g_n) \to (M_n , h_n)$ from Lemma \ref{lem:diffeo}. The result now easily follows by the universal property for direct limits \cite{WO93} since the composition of the top and bottom rows are the connecting maps in the direct limits $\A(M, g, \F)$ and $\A(M, h, \F)$, respectively.
\EPf

The $\cs$-algebra $\A(M, g, \F)$ does indeed depend on the choice of the augmented Fredholm filtration $\F = (M_n, U_n)_{n=k}^\infty$. However, we will see in the next section that the $K$-theory groups of $\A(M, g, \F)$ do not depend on the choice of the tubular neighborhoods $\{U_n\}_{n=k}^\infty$ and, moreover,  if $M$ has an appropriate spin structure then the $K$-theory groups do not depend on the choice of filtrating manifolds $\{M_n\}_{n=k}^\infty$.

We will now consider two examples from the literature that are directly related to this construction.

\begin{example}\label{ex:HKT}
Consider $M = \E$, with metric $g$ induced by the inner product $\< \cdot, \cdot  \>$, and Fredholm filtration given by a flag $\{ E_n\}_n$ of $\E$ as in Example 3.4 (i). Setting $\nu E_n = U_n = E_{n+1}$, we obtain an augmented Fredholm filtration $\F = (E_n, E_{n+1})_{n=1}^\infty$ of $\E$.  We thus have the $\cs$-algebra
\begin{equation*}
\A(\E, g, \F) = \varinjlim\, \A(E_n)
\end{equation*}
as constructed above. Since $TE_n \cong E_n \times E_n$ is trivial, we have that
\begin{equation*}
\A(E_n) \cong \S \gtimes C_0(E_n, \Cliff(E_n)) = \SC(E_n)
\end{equation*}
as in Definition 3.1 of Higson-Kasparov-Trout \cite{HKT98}. Also, since $\nu E_n = U_n = E_{n+1}$, it follows that the connecting map $\alpha_n : \A(E_n) \to \A(E_{n+1})$ can be canonically identified with the Bott periodicity map
\begin{equation*}
\beta_{(n+1)n} = \alpha_n : \A(E_n) \to \A(E_{n+1})
\end{equation*}
of Definition 3.1 in  \cite{HKT98}. 
Using an approximation argument to deal with the dense subalgebra of compactly supported functions, it follows that the $\cs$-algebra $\A(\E, g, \F)$ is isomorphic to the $\cs$-algebra 
\begin{equation*}
\A(\E) = \varinjlim_{E_a \subset \E} \A(E_a)
\end{equation*}
where the direct limit is taken over the directed system of {\it all} finite dimensional subspaces $E_a \subset \E$. See also Lemma 2.6 and the discussion after  Definition 4.6 of Higson-Kasparov \cite{HK01}. This $\cs$-algebra has important applications to the Baum-Connes and Novikov Conjectures \cite{HK01, HKT98, Yu00}
\end{example}

\begin{example}
Another example, which generalizes the previous one, comes from the Thom isomorphism theorem for infinite rank Euclidean vector bundles  \cite{Trou03}. Suppose $M$ is the total space of a smooth (locally trivial)  vector bundle $p : M \to X$, with fiber $\E$ and structure group $GL(\E)$, over a smooth, {\it finite-dimensional} Riemannian manifold $X$ of dimension $k$. Since the fiber $\E$ is infinite-dimensional, we may assume \cite{Dix63} that $M = X \times \E$ is trivial. The inner product $\langle \cdot, \cdot \rangle$ on $\E$ then canonically induces a Euclidean metric structure on the bundle $M$. Using the isomorphism
\begin{equation*}
TM \cong TX \times T\E = TX \times (\E \times \E)
\end{equation*}
we canonically endow the total space $M$ with the structure of a Riemannian Hilbert manifold. Also, since $TM$ is trivial, it follows that $M$ has a canonical structure as a Fredholm manifold.

Let $\{E_n\}_{n=1}^\infty$ be a flag for $\E$. For each $n \geq k+1$, let 
\begin{equation*}
M_n = X \times E_{n-k} \to X
\end{equation*}
denote the trivial vector subbundle of rank $n-k$. One can then check that the collection of submanifolds $\{M_n\}_{n = k+1}^\infty$  determines a Fredholm filtration of $M$ such that we can canonically identify
the total space $\nu M_n$ of the normal bundle of $M_n \hookrightarrow M_{n+1}$ as $M_{n+1}$. We then have that $\F = (M_n, M_{n+1})_{n=k+1}^\infty$ is an augmented Fredholm filtration for $M$. Since $M_n = X \times E_{n-k}$ we have 
$$\A(M_n) \cong \A(E_{n-k}) \gtimes \CC(X) \cong \S \gtimes \CC(E_{n-k}) \gtimes \CC(X).$$ It follows from Proposition 2.13, the results in \cite{Trou03}, and a similar approximation argument that
$$\A(M, g, \F) \cong  \A(\E) \gtimes \CC(X) \cong \A(M, \nabla_0, X)$$
where $\A(M, \nabla_0, X)$ is the $\cs$-algebra of the affine Euclidean bundle $p : M \to X$, equipped with the trivial connection $\nabla_0 = d$, as in Definition 3.11 of  \cite{Trou03}.

\end{example} 

\section{K-theory, Spin Structures and Poincar\'e Duality}

\newcommand{\so}{\mathcal{SO}}
\newcommand{\spin}{\text{{Spin}}}

In this section we discuss the relationship between the topological $K$-theory groups, the (compactly supported) $K$-homology groups of a Fredholm manifold $M$ and the $K$-theory groups of the $\cs$-algebra $\A(M, g, \F)$ we constructed in the last section. When an oriented Riemannian Fredholm manifold $M$ has been equipped with an appropriate infinite-dimensional spin structure, we will see that all of these groups coincide, as in the finite-dimensional spin manifold setting.

\subsection{The topological $K$-theory of a Fredholm manifold.}

Mukherjea \cite[Sec.2]{Mkhr2}, in the context of generalized cohomologies 
obtained from a spectrum on the category of compact spaces, defined the
corresponding cohomology groups for Fredholm manifolds. Based on his work, 
we are led to make the following definition.

\begin{dfn}\label{D:topKth}
Let $M$ be smooth Fredholm manifold with augmented Fredholm filtration $\F = (M_n, U_n)_{n=k}^\infty$.
The {\it $j^{th}$ topological $K$-theory group of $(M, \F)$}, denoted $K^{\infty-j}(M,\F)$,
is defined to be the direct limit
$$
K^{\infty-j}(M,\F)=\varinjlim\,  K^{n-j}(M_n), \text{ for } j=0,1,
$$
where the connecting maps are the Gysin (or shriek) maps \cite{CS84, Kar78}
$$(j_n)_! : K^{n-j}(M_n) \to K^{n+1 - j}(M_{n+1})$$
associated to the inclusions $j_n : M_n \hookrightarrow M_{n+1}$. These may be
obtained from diagram (\ref{eqn:top.dirlim}),
via the functoriality properties of topological $K$-theory, 
as the composition of Gysin maps
\begin{equation}
\xymatrix{ 
K^{n+1-j}(\nu M_n) \ar[r]^{(\phi_n)_{!} }_{{\cong}} & K^{n+1-j}(U_n) \ar[rr]^{(k_n)_{!}} 
                                      & \quad    & K^{n+1-j}(M_{n+1}) \\
K^{n-j}(M_n) \ar[u]^{s_{!}=\text{Thom}}_{\cong} \ar@{-->}[urrr]_{(j_n)_{!}}
}
\end{equation}
where the map $s_!$ is the Gysin map associated to the zero section $s : M_n \to \nu M_n$, and which induces the Thom isomorphism. (Compare this with diagram (\ref{eqn:alg.dirlim}).)
\end{dfn}
Clearly, the definition of the topological $K$-theory of $M$ does not depend on the choice of tubular neighborhoods $\{U_n\}_n$ (or any Riemannian metric $g$) but does, {\it a priori},  depend
on the choice of Fredholm filtration $\{M_n\}_n$, as does the definition of $\A(M, g, \F)$. However, if $M$ has a certain infinite-dimensional spin structure, then these topological $K$-theory groups $K^{\infty -j}(M, \F)$ do not depend on the choice of $\F = (M_n, U_n)_n$.
%

\subsection{Fredholm Spin$_q$-structures}
Recall the notation introduced at the beginning of Section~3. Let $\E$ be a separable infinite-dimensional Euclidean space. For $1 \leq q \leq \infty$, let $GL_q(\E) = GL(\E) \bigcap (I+P_q)$, where $P_q$ is the $q$-th Schatten-von Neumann perturbation class. Let  $\O(\E)$  denote the orthogonal operators on $\E$. We let $\O_q(\E) = \O(\E) \bigcap GL_q(\E)$ and 
let $\S\O_q(\E)$ denote the connected component of $I$ in $\O_q(\E)$. All of these groups are infinite-dimensional Banach-Lie groups \cite{dlH72} with manifold topology given by the restriction of the norm $\| \cdot \|_q$. Note that since $P_q \subset \K(\E)$, it follows that $GL_q(\E) \subset GL_\K(\E)$ and so any Hilbert manifold with $GL_q(\E)$ as structure group has a canonical Fredholm structure as in Definition \ref{def:fredmfld}.

Let $M$ be a smooth, paracompact, connected Hilbert manifold, without boundary, modeled on $\E$. Let $\xi : E \to M$ be a smooth (locally trivial) vector bundle over $M$, with fiber $\E$, endowed with a reduction of the structure group from $GL(\E)$ to $GL_q(\E)$. A {\it Riemannian $q$-structure} \cite[Def 2.1]{Anst} on $\xi$ is a reduction of the structure group from $GL_q(\E)$ to $\O_q(\E)$. Since $M$ is paracompact, this may be accomplished by using a partition of unity to define a smooth bundle metric $g_x$ on the fibers $E_x$ of $\xi$. If $\xi$ is the tangent bundle $\pi : TM \ra M$, with Fredholm structure group $GL_q(\E)$, then we say that $M$ has a Riemannian $q$-structure.

\begin{dfn}\cite[Def 2.2]{Anst}
A Riemannian $q$-structure on $\xi : E \to M$ is \emph{orientable} if $\xi$ admits a further reduction of its structure group to $\S\O_q(\E)$. A given reduction will be called an {\it orientation} and $\xi$ will be said to have an \emph{oriented Riemannian $q$-structure}. 
\end{dfn}

A proof of the following can be found in  \cite[Prop 6.2]{Ksch} or \cite[Thm 2.1]{Anst}.
\begin{thm}
A Riemannian $q$-structure on $\xi : E \ra M$ is orientable if and only if the first Stieffel-Whitney class $w_1(\xi) \in H^1(M, \Z_2)$ vanishes. In particular, if $M$ has a Riemannian $q$-structure, then $M$ is orientable if and only if $w_1(M) = w_1(TM) = 0$.
\end{thm}

For the theory of Stieffel-Whitney classes associated to Hilbert bundles over Hilbert manifolds that we are considering, see Koschorke \cite{Ksch}. Note that, contrary to the finite-dimensional case, these characteristic classes are {\it not} diffeomorphism invariants, in general. (See \cite[Example 6.2]{Ksch} for details.) 

Since $\S\O_q(\E)$ is of index $2$ in $\O_q(\E)$, it follows that the universal covering $\Spin_q(\E)$ is a Banach-Lie group and the covering map is $2$-sheeted. We thus have an exact sequence of (paracompact) topological groups
\begin{equation*}
\xymatrix{
1 \ar[r] &  \Z_2 \ar[r] & \Spin_q(\E) \ar[r]^{\rho} &  \S\O_q(\E) \ar[r] & 1}
\end{equation*}
Concrete realizations of these infinite-dimensional spin groups were constructed for $q=1$ by P.~de~la~Harpe \cite{dlH} and for $q=2$ by Plymen and Streater \cite{PlSt}. However, we will not need explicit constructions of these spin groups, only the fact that they are $2$-sheeted covering groups of the associated special orthogonal groups, as in the finite-dimensional case. In the following, we may abbreviate $\spin_q$ and $\so_q$ for $\spin_q(\E)$ and $\so_q(\E)$, respectively.

\begin{dfn} $($\cite[Def 2.4]{Anst}$)$
Suppose $\xi : E \to M$ has an $\so_q$-structure, \ie an oriented Riemannian $q$-structure.
A \emph{$\spin_q$-structure} on $\xi$  is a principal bundle extension associated to 
the covering map
\begin{equation*}
{\rho} : \spin_q \to \so_q
\end{equation*}
of the principal $\so_q$-bundle of linear frames of $\xi$. If $M$ is a Fredholm manifold with oriented Riemannian $q$-structure, then a \emph{$\spin_q$-structure} on $M$ is a $\spin_q$-structure on $\pi : TM \ra M$. We will then call $M$ a \emph{Fredholm $\spin_q$-manifold}.
\end{dfn}

That is, if $p : L \to M$ is the principal  $\so_q$-bundle of oriented orthonormal frames of $\xi : E \to M$, then a $\spin_q$-structure for $\xi$ is a principal $\spin_q$-bundle $q : \Sigma \to M$ such that $\Sigma$ is a $2$-fold covering of $L$, the restriction of the covering map $\tilde{\rho} : \Sigma \to L$ to the fibers are $2$-sheeted coverings and $$\tilde{\rho}(s \cdot g) = \tilde{\rho}(s) {\rho}(g) \text{ and }  q(s) = p(\tilde{\rho}(s))$$ for all $s \in \Sigma$ and $g \in \spin_q$. Thus, the following diagram commutes:
\begin{equation*}
\xymatrix{
\Sigma \ar[r]^{\tilde{\rho}} \ar[d]^q & L \ar[d]_p \\
M \ar[r]^{\text{Id}_M} & M}
\end{equation*}

For $q =1$ de la Harpe has shown that the existence of a $\spin_q$-structure on a Fredholm manifold $M$ with oriented Riemannian $q$-structure is equivalent to the vanishing $w_2(M) = 0$ of the second Stieffel-Whitney class in $H^2(M, \Z_2)$. We wish to extend his result to all values $1 \leq q \leq \infty$ and all $\so_q$-vector bundles. Although his argument for $q =1$ almost certainly holds in the general case, we will provide a more direct proof using  an argument of Lawson and Michelson \cite{LM89} from the finite-dimensional spin case. In order to do that, we need the following cohomology computation, which follows from some results in the literature \cite{dlH72,dlh2}, but we provide a proof for completeness.

\begin{lem} 
For $1 \leq q \leq \infty$, $H^1(\S\O_q(\E), \Z_2) \cong \Z_2.$
\end{lem}

\Pf Choose a flag $\{E_n\}$ for $\E$ as in Example 3.4 (i). This induces an inclusion of topological groups
\begin{equation*}
SO(\infty) = \varinjlim\,  SO(n) \hookrightarrow \S\O_q(\E)
\end{equation*}
which, by Proposition 3 in  \cite{dlh2}, is a homotopy equivalence. Hence, using the identity as basepoint, we have by  Bott periodicity \cite{Bot59}:
\begin{equation*}
\pi_1(\S\O_q(\E)) \cong \pi_1(SO(\infty)) \cong  \varinjlim\, \pi_1(SO(n)) \cong \Z_2.
\end{equation*}
Since $\S\O_q(\E)$ is connected with abelian fundamental group, it follows that
\begin{equation*}
H_1(\S\O_q(\E)) \cong \pi_1(\S\O_q(\E)) \cong \Z_2.
\end{equation*}
The result now follows from the Universal Coefficient Theorem in cohomology:
\begin{equation*}
H^1(\S\O_q(\E), \Z_2) \cong \Hom(H_1(\S\O_q(\E), \Z), \Z_2) \cong \Hom(\Z_2, \Z_2) \cong \Z_2
\end{equation*}
and we are done. \hfill \EPf

\begin{thm}\label{thm:spinchar}
Let $\xi : E \to M$ be a Hilbert bundle with oriented Riemannian $q$-structure. Then $\xi$ has a $\spin_q$-structure if and only if the second Stieffel-Whitney class $w_2(\xi) \in H^2(M, \Z_2)$ vanishes. In particular, if $M$ is a Fredholm manifold with oriented Riemannian $q$-structure, then there exists a $\spin_q$-structure on $M$ if and only if $w_2(M) = 0$.
\end{thm}

For the following, recall that in principal bundle theory, if $M$ is a paracompact space and $G$ is a topological group, then $H^1(M, G)$ is isomorphic to the set of isomorphism classes of principal $G$-bundles on $M$, where we are using $\check{\text{C}}$ech cohomology. (See Appendix A of Lawson and Michelsohn \cite{LM89}.)

\Pf Let $p: L \to M$ be the principal $\so_q$-bundle of oriented orthonormal frames of $\xi$. We then have a fibration
\begin{equation*}
\xymatrix{
\S\O_q(\E) \ar[r]^-{i} & L \ar[r]^-{p} & M}
\end{equation*}
which induces an exact sequence
\begin{equation*}
\xymatrix{
H^1(M, \Z_2) \ar[r]^-{p^*} & H^1(L, \Z_2) \ar[r]^-{i^*} & H^1(\S\O_q(\E), \Z_2) \ar[r]^-{\delta_\xi} & H^2(M, \Z_2)}
\end{equation*}
in $\check{\text{C}}$ech cohomology. It follows by the above discussion (see also  \cite[Thm 2.3]{Anst}) that $\xi$ has a $\spin_q$-structure if and only if there is a cohomology class $\alpha = \alpha(\xi) \in H^1(L, \Z_2)$ such that $i^*(\alpha) \neq 0$ since a $\spin_q$-structure on $\xi$ determines a nontrivial $2$-sheeted covering of $L$. Let $g_2$ be the generator of $H^1(\S\O_q(\E), \Z_2) \cong \Z_2.$ It follows that  $\xi$ has a $\spin_q$-structure if and only if there is a cohomology class $\alpha = \alpha(\xi) \in H^1(L, \Z_2)$ such that $i^*(\alpha) = g_2$. Consequently,  by exactness of the sequence above, we have that this holds if and only if
\begin{equation*}
w_2(\xi) = \delta_\xi(g_2) = \delta_\xi(i^*(\alpha)) = 0 \in H^2(M, \Z_2).
\end{equation*}
The fact that the second Stieffel-Whitney class of $\xi$ is given by 
\begin{equation*}
w_2(\xi) = \delta_\xi(g_2) \in H^2(M, \Z_2)
\end{equation*}
 follows from the universal properties of these classes \cite[Proposition 6.3]{Ksch}. \hfill \EPf

Consequently, if $\xi : E \to M$ admits a $\spin_q$-structure determined by $\alpha(\xi) \in H^1(L, \Z_2)$ then the most general $\spin_q$-structure on $\xi$ is of the form $\alpha(\xi) + p^*(\beta)$ where $\beta \in H^1(M, \Z_2)$. Thus, there is a bijection between the set of (isomorphism classes of) $\spin_q$-structures on $\xi$ and $H^1(M, \Z_2)$. It follows that a $\spin_q$-structure on $\xi$ (or $M$) is unique if $H^1(M, \Z_2) = 0$.

The next two results are immediate corollaries (see Theorems 2.5 and 2.6  of  \cite{Anst}.)

\begin{prop}\label{T:2out3}
Given $\spin_q$-structures on two out of the three vector bundles 
$\xi_1$, $\xi_2$, and $\xi_1\oplus \xi_2$ on $M$, there is 
a uniquely determined $\spin_q$-structure on the third.
\end{prop}

\begin{prop}\label{T:indSpin}
If $\xi : E \ra M$ admits a $\spin_q$-structure and $f : N \ra M$ is smooth,
then  the pull-back vector bundle $f^*\xi : f^*E \ra N$ admits a $\spin_q$-structure.
\end{prop}

In the context of Fredholm manifolds, the above give:

\begin{cor}\label{cor:FFspin}
Let $M$ be a Fredholm $\spin_q$-manifold. If  $\{ M_n \}_n$ is any associated Fredholm filtration of $M$  then each $M_n$ has a canonical $($finite-dimensional$)$ spin structure.
\end{cor}

\noindent
Indeed, associated to the inclusion $i_n : M_n \ra M$ we have a split short exact sequence
$$
0 \ra TM_n \ra TM|_{M_n} \ra \mu M_n \ra 0.
$$
The normal bundle $\mu M_n$ has a $\spin_q$-structure being trivial, and 
$TM|_{M_n} = i_n^*(TM)$ has one because of Proposition \ref{T:indSpin}. Thus, we have
\begin{equation*}
w_2(\mu M_n) = w_2(TM|_{M_n}) = 0
\end{equation*}
and finally Proposition \ref{T:2out3} gives the result since $w_2(M_n) = 0$.

We end this subsection about spin structures with an example coming 
from certain based loop groups.

\begin{example}
Consider a compact, connected, simply connected, simple Lie group $G$.
Let $\Omega_s G= H^s_0(S^1,G)$ be the
{\it group of based loops on $G$}, \ie  maps from the circle
to $G$ in the $s$th Sobolev space $H^s$ which take a fixed point on $S^1$ into the identity element of 
$G$, where $s \geq 1/2$.
$\Omega_s G$ is a (real) Hilbert Lie group.
%

D. Freed constructed in \cite[Sec.5]{Freed88} a particular Fredholm 
1-structure, coming from a classifying map
$$
\Omega_s G \rightarrow BGL(\infty; \C) \sim \Phi_0,
$$
where $\Phi_0$ denotes the Fredholm operators of index zero.
The resulting frame bundle was called the
{\it geometric frame bundle}. He concluded that the realification of this
geometric frame bundle is trivial and that the Stieffel-Whitney classes 
of $\Omega_s G$
vanish (\cite[Thm 5.30]{Freed88}).
Our Theorem \ref{thm:spinchar} now shows that this is the unique $\spin_1$-structure 
on $\Omega_s G$.
Indeed, the hypothesis on $G$ implies that
$\pi_0(G)=\pi_1(G)=\pi_2(G)=0$, and $\pi_3(G)=\Z$. Consequently
$H_1(\Omega_s G,\Z)=0$ and $H_2(\Omega_s G,\Z) \cong \pi_2(\Omega_s G) \cong 
\pi_3(G) = \Z$.
These imply that
$H^1(\Omega_s G,\Z_2)=0$ and $H^2(\Omega_s G,\Z_2) = \Z/2$.
As $w_2(\Omega_s G)=0$ by Freed's Corollary 5.31, and as 
$\spin_1$-structures
on $\Omega G$ are parametrized by $H^1(\Omega_s G,\Z_2) = 0$, we obtain the 
claimed
uniqueness of the $\spin_1$-structure on $\Omega_s G$.
Moreover, Freed's Fredholm structure is actually the unique
$\spin_q$-structure, for all $1 \leq q \leq \infty$.
\end{example}


\subsection{$K$-homology and Poincar\'e duality}
Recall that if $X$ is a compact space then the $j$-th $K$-homology group of $X$ is the abelian group $K_j(X) = KK^j(C(X), \C)$ which is dual to the $j$-th $K$-theory group $K^j(X) \cong KK^j(\C, C(X))$. The map $X \mapsto K_j(X)$ defines a generalized homology theory on the category of compact spaces and continuous maps \cite{BDF77,Kas81,HR00}.

\begin{dfn}\label{def:K-hmlgy}
Let $M$ be a paracompact space. The \emph{ $j^{th}$ compactly supported $K$-homology group of $M$} is
$$
K_j^c(M)= \varinjlim_{X\subset M} K_j(X),
$$
where the direct limit is over all the {\it compact} subsets $X\subset M$,
and $j=0,1$.
\end{dfn}
\noindent 

In order to prove our Poincar\'e duality result, we need the following result, whose proof requires the $KK$-theory for pro-$\cs$-algebras developed by Weidner \cite{Weid89} and Phillips \cite{Phi89-2}. A heuristic proof would be that since $M \sim M_\infty = \varinjlim \, M_n$, we have in compactly supported $K$-homology that $K^c_j(M) \cong K^c_j(M_\infty) \cong \varinjlim \, K^c_j(M_n)$.

\begin{prop}\label{prop:K-hmlgy}
Let $M$ be a smooth Fredholm manifold. If $\{M_n\}_{n=k}^\infty$ is any Fredholm filtration of $M$ then there is an isomorphism of abelian groups
\begin{equation*}
K_j^c(M) \cong \varinjlim\, K_j^c(M_n), \quad j = 0, 1,
\end{equation*}
where the connecting map $K_j^c(M_n) \to K_j^c(M_{n+1})$ in the direct limit is induced by the inclusion $M_n \hookrightarrow M_{n+1}$.
\end{prop}

\Pf Let $g$ be a Riemannian metric on $M$ compatible with the topology (which exists via paracompactness). Thus, $(M, g)$ is a metric space. Since metric spaces are compactly generated \cite[I.4.3]{White78}, it follows that the algebra $C(M)$ of all continuous complex-valued functions on $M$, with the topology of uniform convergence on compact subsets, is a pro-$\cs$-algebra with involution given by pointwise complex conjugation \cite[Ex 1.3.3]{Phi88}. Let $\CC_M$ denote the collection of all compact subsets $X$ of $M$ ordered by inclusion. Since $M$ is regular, it is completely Hausdorff \cite[Def 2.2]{Phi88}, and so by Corollary 2.9 of \cite{Phi88}, it follows that there is an isomorphism
\begin{equation}\label{eqn:M}
C(M) \cong \varprojlim_{X \in \CC_M} \, C(X)
\end{equation}
of pro-$\cs$-algebras. Similarly, for each $n$, we have an isomorphism
\begin{equation}\label{eqn:M_n}
C(M_n) \cong \varprojlim_{K_n \in \CC_{M_n}} \, C(K_n)
\end{equation}
of pro-$\cs$-algebras where $\CC_{M_n}$ denotes the set of all compact subsets $K_n$ of $M_n$ ordered by inclusion. Let $M_\infty = \bigcup_n M_n = \varinjlim \, M_n$ with the direct limit topology. Since $M_\infty$ is countably compactly generated in the direct limit topology, we then have an isomorphism
\begin{equation}\label{eqn:Minfty}
C(M_\infty) \cong \varprojlim_{n} C(M_n)
\end{equation}
of pro-$\cs$-algebras. By Theorem \ref{Thm:FredFilt} the inclusion $M_\infty \hookrightarrow M$ is a homotopy equivalence, hence the pro-$\cs$-algebras $C(M)$ and $C(M_\infty)$ have the same homotopy type.

Using the fact that Weidner's $KK$-groups $KK^j_W(A, B)$ for pro-$\cs$-algebras \cite{Weid89,Phi89-2} extend Kasparov's $KK$-groups for $\cs$-algebras \cite{Kas81}, are homotopy-invariant, and convert inverse limits to direct limits \footnote{Note that there is a typo in the statement of  \cite[Thm 5.1]{Weid89}} in the $K$-homology variable, we compute as follows:
\begin{equation*}
\begin{aligned}
K_j^c(M) & = \varinjlim_{X \in \CC_M} KK^j(C(X), \C) & & \text{(Definition \ref{def:K-hmlgy})}   \\
& \cong KK^j_W(\varprojlim_{X \in \CC_X} C(X), \C)  && \text{(By \cite[Thm 5.1]{Weid89})} \\
& \cong  KK^j_W(C(M), \C) && \text{(By  Eqn (\ref{eqn:M}))}  \\
& \cong KK^j_W(C(M_\infty), \C)  && \text{(homotopy invariance)} \\
& \cong  KK^j_W(\varprojlim_{n} \, C(M_n), \C) &&  \text{(By  Eqn (\ref{eqn:Minfty}))}  \\
& \cong \varinjlim_{n} \, KK^j_W(C(M_n), \C) && \text{(By \cite[Thm 5.1]{Weid89})}  \\
& \cong \varinjlim_{n} \, KK^j_W(\varprojlim_{K_n} \, C(K_n), \C) & & \text{(By  Eqn (\ref{eqn:M_n}))}  \\
& \cong \varinjlim_{n} \, \varinjlim_{K_n} \, KK^j(C(K_n), \C) &&   \text{(By \cite[Thm 5.1]{Weid89})} \\
& \cong \varinjlim_{n} \, K^c_j(M_n) &&  \text{(Definition \ref{def:K-hmlgy})} 
\end{aligned} 
\end{equation*} \EPf

Compare the following result for Fredholm $\spin_q$-manifolds with \cite[Thm 2.1]{Mkhr2}.

\begin{thm}[{\bf Poincar\'e duality}]\label{Thm:FMPD}
If $M$ is a smooth Fredholm $\spin_q$-manifold with augmented Fredholm filtration $\F$, there is an isomorphism
$$
K^{\infty-j}(M,\F) \cong K_j^c(M)  
$$
\end{thm}

\Pf Let $\F = (M_n, U_n)_{n=k}^\infty$ be the augmented Fredholm filtration. Since $M$ is a Fredholm $\spin_q$-manifold, each $M_n$ has a canonical spin structure by Corollary \ref{cor:FFspin}. By \cite[Cor 31]{Ren03} (or \cite[Exercise 11.8.11]{HR00}) we have a natural Poincar\'e duality isomorphism
\begin{equation*}
\xymatrix{P_n : K^{n-j}(M_n) \ar[r]^-{\cong} &  K^c_j(M_n)}
\end{equation*}
given by the cap product with the fundamental class $[M_n]$. Naturality is the assertion that the Poincar\'e duality diagram
\begin{equation*}
\xymatrix{K^{n-j}(M_n) \ar[r]^-{{j_n}_!} & K^{n+1-j}(M_{n+1}) \\
K_j^c(M_n) \ar[u]-^{P_n}_{\cong} \ar[r]^-{{j_n}_*} & K_{j}^c(M_{n+1}) \ar[u]^{\cong}_{P_{n+1}}}
\end{equation*}
commutes, where $j_n : M_n \hookrightarrow M_{n+1}$. It now follows that:
$$
\begin{aligned}
K_j^c(M)
&\cong \varinjlim\,  K_j^c(M_{n}) && \text{(Proposition \ref{prop:K-hmlgy})} \\
&\cong \varinjlim\,  K^{n-j}(M_{n}) && \text{(classical Poincar\'e duality)} \\
&= K^{\infty-j}(M,\F) & & \text{(Definition \ref{D:topKth})}
\end{aligned}
$$ 
as desired. \EPf

\subsection{$K$-theory of the $\cs$-algebra $\A(M, g, \F)$}

First we discuss the finite dimensional results we will need. Let $M_n$ be an oriented Riemannian $n$-manifold.
An important relationship between the non-commutative $\cs$-algebra
$\CC(M_n) = C_0(M_n, \Cliff (TM_n))$ and the commutative $\cs$-algebra $C_0(M_n)$ is given by
\spinc-structures~\cite{LM89}. Let $\C_1 = \Cliff(\R)$ denote the first complex
Clifford algebra. The following is adapted from
Theorem 2.11 of Plymen \cite{Ply86} and Proposition II.A.9 of Connes
\cite{Con94}.

\begin{prop}\label{prop:morita}
If $n = 2k$ is even, there is a bijective correspondence between
\spinc-structures on $M_n$ and Morita equivalences $($in the sense of Rieffel
\cite{RW98,Rie82}$)$ between the $\cs$-algebras $C_0(M_n)$ and
$\CC(M_n)$. Thus, $\A(M_n)$ is Morita equivalent to $C_0(\R \times M_n)$. If $n
= 2k+1$ is odd, then
\spinc-structures on $M$ are in bijective correspondence with Morita
equivalences $C_0(M_n) \sim \CC(M_n) \gtimes \C_1$.
\end{prop}

Although $\CC(M_n)$ and $\A(M_n)$ carry natural $\Z_2$-gradings, when we
consider their $\cs$-algebra $K$-theory, we will ignore these
gradings. That is, if $A$ is any $\cs$-algebra --- graded or not --- then $K_j(A)$
($j = 0, 1$) will denote the $K$-theory group of the underlying $\cs$-algebra,
without the grading. Since $\cs$-algebra $K$-theory is Morita invariant, we have the following.

\begin{cor}\label{cor:topKthry}
If $M_{2k}$ is an even-dimensional oriented Riemannian manifold with
\spinc-structure, there is a canonical $K$-theory 
isomorphism\footnote{It is also true that $K_j(\CC(M_{2k})) \cong K^j(M_{2k})$, 
but we shall not use this here.}
\begin{equation*}
K_j(\A(M_{2k})) \cong K^{j+1}(M_{2k}).
\end{equation*}
\end{cor}

The next result is proved by Trout  \cite[Thm 2.14]{Trou03}:
\begin{ThomThm}\label{thm:Thomfrc}
If $E \to M_n$ is a smooth finite-rank affine Euclidean bundle, then the
$*$-homomorphism $\Psi_p : \A(M_n) \to \A(E)$ from Theorem \ref{thm:Thomhom} induces an isomorphism of abelian
groups:
$$\Psi_{*} : K_j(\A(M_n)) \to K_j(\A(E)), \text{ for } j = 0, 1.$$
\end{ThomThm}

\noindent
In fact, it is the $\cs$-algebraic formulation of the classical Thom isomorphism $\Phi : K^j(M) \to K^j(E)$ from topological $K$-theory.

\begin{cor}\cite[Cor 2.20]{Trou03}
If $E$ is a finite even-rank oriented Euclidean \spinc-bundle $($with spin
connection $\nabla$$)$ on an even-dimensional
oriented Riemannian \spinc-manifold $M_n$, then $\Psi_p : \A(M_n) \to \A(E)$
induces the topological Thom isomorphism $\Phi$, as depicted in the 
following commutative diagram:
$$\xymatrix{
K_j(\A(M_n)) \ar[r]^{\Psi_*} \ar[d]_{\cong} & K_j(\A(E)) \ar[d]^{\cong} \\
K^{j+1}(M_n) \ar[r]^{\Phi} & K^{j+1}(E)
}$$
\end{cor}

Although the connecting maps $\alpha_n : \A(M_n) \to \A(M_{n+1})$ are not functorial at the $\cs$-algebra level (as in diagram (\ref{diag:functor})), they are at the level of $K$-theory.

\begin{lem} The following diagram of abelian groups
\begin{equation}\label{diag:Kthryfunc}
\xymatrix{ & K_j(\A(M_{n+1})) \ar[dr]^-{(\alpha_{n+1})_*} &  \\ K_j(\A(M_n)) \ar[ur]^-{(\alpha_{n})_*} \ar[rr]_-{(\alpha_n^{n+2})_*}  &  & K_j(\A(M_{n+2}))}
\end{equation}
commutes for all $n \geq k$ and $j = 0, 1$, where $\alpha_n^{n+2} : \A(M_n) \to \A(M_{n+2})$ is any Gysin map induced by the inclusion $M_n \hookrightarrow M_{n+2}$ $($as in Diagram $(6))$.
\end{lem}

\Pf The functor $M_n \mapsto K_j(\A(M_n))$ from the category of finite-dimensional smooth (Riemannian) manifolds is homotopy-invariant, has Gysin maps (independent of the choice of tubular neighborhood) and, most importantly, a {\it transitive} Thom homomorphism \cite[Lem 3.10]{Trou03}. The result now follows from the corresponding proof in Karoubi \cite[Props 5.22 and 5.24]{Kar78} for topological $K$-theory. \EPf

We now come to the main result of our paper.

\begin{thm}\label{T:Kth-A(M)}
Let $M$ be a smooth Fredholm $\spin_q$-manifold with Riemannian metric $g$ and augmented Fredholm filtration $\F = (M_n, U_n)_{n=k}^\infty$. With a dimension shift,
the $K$-theory of $\A(M,g,\F)$ coincides with the
topological $K$-theory of $(M, \F)$ and the (compactly supported) $K$-homology of $M$:
\begin{equation*}
K_{j+1}(\A(M, g, \F)) \cong K^{\infty-j}(M, \F) \cong K_j^c(M).
\end{equation*}
\end{thm}

\Pf Indeed, using the fact that $2 \Z$ is cofinal in $\Z$, we can restrict to the
even-dimensional subsequences in the directed limits under
consideration:
\begin{equation*}\label{eqn:KthM}
\begin{aligned}
K^{\infty-j}(M,\F) 
&=     \varinjlim\,  K^{n-j}(M_n) \qquad\qquad && \text{(Definition \ref{D:topKth})} \\
&\cong \varinjlim\,  K^{2n-j}(M_{2n}) && \text{(cofinal property of direct limits)} \\
&\cong \varinjlim\,  K^{j+2}(M_{2n}) && \text{(Bott periodicity)} \\
&\cong \varinjlim\,  K_{j+1}(\A(M_{2n})) && \text{(Corollary \ref{cor:topKthry})} \\
&\cong \varinjlim\,  K_{j+1}(\A(M_{n})) &&\text{(cofinal property of direct limits)} \\
&\cong K_{j+1}(\varinjlim\,  \A(M_n)) && \text{(continuity of $K$-theory)}   \\
& = K_{j+1}(\A(M,g,\F)) && \text{(Definition \ref{def:C*algFM})}
\end{aligned}
\end{equation*}
\EPf

\noindent
As the compactly supported $K$-homology of $M$ does not depend on
the metric and on the choice of augmented filtration,  
we get in particular the following independence on the metric and the filtration
(compare again with \cite[Thm 2.1]{Mkhr2}):
\begin{cor}
If $M$ is a smooth Fredholm $\spin_q$-manifold, as above, 
then its topological $K$-theory $K^{\infty-j}(M, \F)$ and the $K$-theory of $\A(M,g,\F)$
do not depend on the choices of the metric $g$ and augmented Fredholm filtration $\F$.
\end{cor}

Another easy consequence is:
\begin{cor} If $\E$ is a separable infinite-dimensional Euclidean space, then
\begin{equation*}
K_j(\A(\E)) \cong K_j(\A(S_{\E})) \cong
\left\{ \begin{array}{r@{ , \quad}l} 
            0  & \text{ if } j=0, \\
            \Z & \text{ if } j=1  \end{array} \right.
\end{equation*}
where $S_\E$ denotes the unit sphere in $\E$.
\end{cor}

\bibliographystyle{amsplain}
\providecommand{\bysame}{\leavevmode\hbox to3em{\hrulefill}\thinspace}
\providecommand{\MR}{\relax\ifhmode\unskip\space\fi MR }
\providecommand{\MRhref}[2]{%
  \href{http://www.ams.org/mathscinet-getitem?mr=#1}{#2}
}
\providecommand{\href}[2]{#2}

\end{document}